\documentclass[reqno,12pt]{amsart}
\usepackage{amsmath, amssymb, amsthm, amsfonts} 
\usepackage[english]{babel}
\usepackage{bbm}
\usepackage{graphicx}
\usepackage{url}
\usepackage{soul}
\usepackage{epstopdf}
\usepackage[ruled,vlined]{algorithm2e}
\usepackage[a4paper,bindingoffset=0.5cm,left=2cm,right=2cm,top=2.5cm,bottom=2cm,footskip=.8cm]{geometry}
\usepackage{rotating}
\usepackage{amsbsy,enumerate}
\usepackage{comment}
\usepackage{mathrsfs} 

\newcommand{\bs}{\boldsymbol}

\newcommand{\vb}{\vspace{3.2mm}}
\renewcommand{\hat}{\widehat}

\newcommand{\diag}{\mathrm{diag}}

\newcommand{\vertiii}[1]{{\left\vert\kern-0.25ex\left\vert\kern-0.25ex\left\vert #1 \right\vert\kern-0.25ex\right\vert\kern-0.25ex\right\vert}}
\allowdisplaybreaks
\usepackage{hyperref}
\newcommand{\e}{\mathbb{E}}
\newcommand{\p}{\mathbb{P}}
\newcommand{\QQ}{\mathbb{Q}}
\newcommand*\diff{\mathop{}\!\mathrm{d}}

\allowdisplaybreaks

\newtheorem{lemma}{Lemma}
\newtheorem{corollary}{Corollary}

\newtheorem{theorem}{Theorem}
\newtheorem{remark}{Remark}

\newtheorem{proposition}{Proposition}

\usepackage[utf8]{inputenc}
\usepackage{type1cm}         
\usepackage{multicol}        
\usepackage[bottom]{footmisc}

\usepackage{newtxtext}       %
\usepackage{newtxmath}       

\usepackage{tikz}
\usetikzlibrary{plotmarks}
\usetikzlibrary{decorations.pathreplacing}
\usepackage{pgfplots}
\pgfplotsset{width=10cm,compat=1.9}
\usepgfplotslibrary{external}
\usepackage{subcaption}
\usepackage{standalone}
\usetikzlibrary{automata,positioning}
\usetikzlibrary{decorations.pathreplacing,decorations.markings,shapes.geometric}
\usetikzlibrary{shapes,arrows}
\usetikzlibrary{backgrounds,calc,positioning}

\begin{document}

	\title[Bankruptcy probabilities under non-Poisson inspection]{Bankruptcy probabilities under non-Poisson inspection}
\author{Florine Kuipers, Michel Mandjes and Sara Morcy}
	
	\begin{abstract}
		This paper concerns an insurance firm's surplus process observed at renewal inspection times, with a focus on assessing the probability of the surplus level dropping below zero. For various types of inter-inspection time distributions, an explicit expression for the corresponding transform is given. In addition, Cram\'er-Lundberg type asymptotics are established. Also, an importance sampling based Monte Carlo algorithm is proposed, and is shown to be logarithmically efficient.

\vb


\vb

\noindent
{\sc Affiliations.} The authors  are with the Korteweg-de Vries Institute for Mathematics, University of Amsterdam, Science Park 904, 1098 XH Amsterdam, The Netherlands. MM is also with the Mathematical Institute, P.O. Box 9512,
2300 RA Leiden,
The Netherlands; E{\sc urandom}, Eindhoven University of Technology, Eindhoven, The Netherlands; Amsterdam Business School, Faculty of Economics and Business, University of Amsterdam, Amsterdam, The Netherlands. MM's research has been funded by the NWO Gravitation project N{\sc etworks}, grant number 024.002.003. 

\noindent Date: {\it \today}.

\vb

\noindent
{\sc Acknowledgments.} 
The second author thanks Onno Boxma (Eindhoven University of Technology) and Offer Kella (Hebrew University, Jerusalem) for helpful discussions, and Hansj\"org Albrecher (University of Lausanne) for useful feedback on the predecessor paper \cite{BM2}.

\vb

\noindent
{\sc Email.} \url{kuipers.florine@gmail.com}, \url{m.r.h.mandjes@uva.nl}, \url{sara0morcy@gmail.com}.

	\end{abstract}

	\maketitle

\section{Introduction}
Ruin theory constitutes a prominent subdomain of actuarial science.  A key model, usually referred to as the {\it Cram\'er-Lundberg model}, describes the surplus process of an insurance company that receives income (at a constant rate) via premiums, and that pays out claims to its clients according to a compound Poisson process. The key quantity studied is the {\it ruin probability}, being the probability that the surplus level drops below~0 (`ruin'), either before a given point in time or over an infinite horizon.
Ruin theory is a mature and active branch within applied probability: for a broad variety of types of surplus processes, covering the above-mentioned Cram\'er-Lundberg model, one has succeeded to evaluate (or approximate) the ruin probability. 
These findings are of great practical interest, as they facilitate the selection of an initial surplus level that renders the ruin probability below a given threshold. 

Mathematically, ruin is usually defined in terms of a certain continuous-time stochastic process reaching a given set (e.g., the surplus level dropping below~0). In many practical situations, however, the stochastic process under considerations is not, or can not, be continuously observed, but is rather inspected at discrete epochs. This motivates the analysis of setups in which the firm goes bankrupt when the surplus level drops below zero {\it at an inspection time}.  

The most basic setup corresponds to Poisson inspections, i.e., there are exponentially distributed `inter-inspection times' at which it is verified whether the surplus level has become negative. This case, with the surplus level evolving according to a Cram\'er-Lundberg model or more general processes, has been dealt with in great detail before \cite{ABEK,BEJ,BM2}. Let $p(u,t)$ denote the probability of ruin at an inspection epoch before time $t$, which in the sequel we refer to as the {\it bankruptcy probability}. Given that the initial surplus level is $u>0$, a main finding concerns an explicit expression for the transform
\[\pi(\alpha,\beta)\coloneqq\int_0^\infty e^{-\alpha u} p(u,T_\beta)\,{\rm d}u,\]
where $T_\beta$ is an exponentially distributed quantity with rate $\beta$, sampled independently of the surplus process. 
The time-dependent bankruptcy probability $p(u,t)$ can be numerically approximated by applying  well-established numerical inversion techniques, such as the fast and accurate procedures presented in  \cite{AW,dI}. In addition,  the asymptotics of the all-time bankruptcy probability $p(u)\coloneqq p(u,\infty)$ (as the initial surplus $u$ grows large) have been identified; see e.g.\ the  results in \cite{BM2} for the case of both light-tailed and heavy-tailed claim sizes. 

Note that, from a practical standpoint, a setup with Poisson inspections is not particularly realistic. The ambition of the present paper is therefore to extend the results for exponentially distributed inter-inspection times to a considerably richer class. 
Acknowledging that the family of {\it phase-type distributions} offers a flexible and conceptually appealing way to approximate any distribution on the positive half-line \cite[Theorem III.4.2]{ASM2}, our focus lies on inter-inspection times stemming from convenient and practically relevant members of this family. 
The net cumulative claim process $Y(\cdot)$, equalling the initial surplus $u$ decreased by the surplus process $X_u(\cdot)$, is assumed to be described by a spectrally-positive L\'evy process such that $Y(0)=0$.

Our objective is to evaluate, for such inter-inspection times, two quantities: (i)~the transform $\pi(\alpha,\beta)$ of the bankruptcy probability, uniquely characterizing the time-dependent bankruptcy probability for any initial surplus $u$, and (ii)~the Cram\'er-Lundberg asymptotics of the all-time bankruptcy probability $p(u)$, characterizing its behavior in the regime that the initial surplus $u$ grows large. 
The types of phase-type distributions considered are mixtures of exponential and Erlang distributions, through which any distribution can (in principle) be approximated arbitrarily accurately. 
Although already covered by e.g.\  \cite{BM2}, we start our analysis by the case of Poisson inspections, to showcase the approach followed, thus increasing the transparency of the derivations for the other, more involved, cases. We also present an importance sampling based Monte Carlo approach to estimate $p(u)$, which can handle in principle any inter-inspection time distribution, and which we show to be logarithmically efficient. The experiments show that replacing the inter-inspection times by a two-moments-fit based, phase-type counterpart yields almost exact results.

We proceed by providing a brief, non-exhaustive overview of related literature. The evaluation of ruin probabilities in the conventional, permanently inspected, Cram\'er-Lundberg context is discussed in great detail in various textbooks, such as  \cite{AA,BMbook}.
The use of our bankruptcy concept was advocated in \cite{AGS}, acknowledging the fact that the insurance firm may continue doing business for a while, even when it is technically ruined.
In \cite{AL}  the bankruptcy probability has been evaluated for the Cram\'er-Lundberg setting with Poisson inspections and exponentially distributed claim sizes, while the later papers~\cite{BEJ,BM2} also cover generally distributed claim sizes. In both \cite{AL} and \cite{BEJ}, the inspection rate was allowed to depend on the current surplus level.
It is noted that in~\cite{BEJ}, and in \cite{ABEK} as well, a related queueing (or inventory) model was analyzed as well, with the special feature that the server works even when there are no customers (or orders), building up storage that is removed at the Poisson inspection epochs. The analysis presented in \cite{BM2} features a clean decomposition, which provides an explicit relation between the model under Poisson inspection and its counterpart under permanent inspection, which has been further generalized in \cite{BKM}. 

Our paper has been organized as follows. Section \ref{MOP} formally describes our model, introduces the objects of study, and recalls a number of useful existing results (such as the Wiener-Hopf decomposition for spectrally-positive L\'evy processes). For exponentially and hyperexponentially distributed inter-inspection times the transform $\pi(\alpha,\beta)$ is computed in Section \ref{Exp_HypExp}, whereas Section~\ref{Erl_HypErl} focuses on Erlang and hyper-Erlang inter-inspection times. In Section \ref{refl} it is pointed out what complications arise when one would try to compute $\pi(\alpha,\beta)$ for more general inter-inspection time distributions. Section~\ref{CLA} identifies the Cram\'er-Lundberg-type asymptotics of $p(u)$ for the types of inter-inspection times considered. Then, in Section \ref{RES}, a provably optimal Monte Carlo method for estimating $p(u)$ is developed and numerical experiments are performed. Section \ref{DC} provides a discussion and concluding remarks. 

\section{Model, objective, and preliminaries}\label{MOP}
We denote by $Y(\cdot)$ a {\it spectrally-positive L\'evy process} with Laplace exponent
\[\varphi(\alpha)\coloneqq  \log {\mathbb E}\,e^{-\alpha Y(1)}\]
and $\psi(\cdot)$ its right inverse. To rule out trivial cases, it is assumed that $Y(\cdot)$ is not a subordinator (i.e., a monotonically increasing or decreasing process). 
Throughout this paper we denote by $\bar Y(\cdot)$ and $\underline Y(\cdot)$ the running maximum and running minimum process, respectively, associated with $Y(\cdot)$. In the context of ruin theory, the process $Y(\cdot)$ represents the {\it net cumulative claim process}, which is equivalent to the initial surplus $u>0$ decreased by the surplus process $X_u(\cdot)$. 

Unlike in the conventional ruin-theoretic setting, in this paper we do not observe the process $Y(\cdot)$  permanently, but rather only according to a renewal inspection process. The inter-inspection times are denoted by the sequence $(\Omega_n)_{n\in{\mathbb N}}$ of independent and identically distributed (i.i.d.) random variables, with the $\Omega_n$ distributed as the generic non-negative random variable $\Omega$, and the corresponding inspection times by
\[\bar\Omega_n\coloneqq \sum_{m=1}^n \Omega_m.\]
The sequence $(\Omega_n)_{n\in{\mathbb N}}$ is assumed to be sampled independently from the net cumulative claim process $Y(\cdot)$.

The goal of this paper is to analyze the (finite and infinite-time) bankruptcy probability at inspection epochs, assuming that the inter-inspection times are distributed as the random variable $\Omega$. The central object of study is
\[p(u,t\,|\,\Omega)\coloneqq  {\mathbb P}\left(\exists n\in{\mathbb N}: \bar\Omega_n\leqslant t, Y(\bar\Omega_n)> u\right),\]
i.e., the probability of bankruptcy before time horizon $t$ only considering the value of $Y(\cdot)$ at inspection epochs; for clarity we have included $\Omega$ in the notation, to stress the dependence on the inter-inspection times. In this paper we consider the case of an exponentially distributed time horizon, from which we can obtain its deterministic horizon counterpart by Laplace inversion. We write $T_\beta$ to denote an exponentially distributed random variable with parameter $\beta$ (i.e., with mean $1/\beta$). By sending $\beta$ to 0 we obtain the infinite-time bankruptcy probability. The random variable $T_\beta$, sampled independently from the net cumulative claim process $\bar Y(\cdot)$ and the inter-inspection times $(\Omega_n)_{n\in{\mathbb N}}$, is often referred to as the {\it killing time}. 

We will intensively use the {\it Wiener-Hopf decomposition} for spectrally-positive L\'evy processes \cite[Section 6.5.2]{KYP}. This states that (i)~$\bar Y(T_\beta)$ and $\bar Y(T_\beta)- Y(T_\beta)$ are independent, (ii)~that for any $\alpha\geqslant 0$ and $\beta>0$,
\begin{equation} \label{xidef}\xi(\alpha,\beta)\coloneqq {\mathbb E}\,e^{-\alpha \bar Y(T_\beta)}=\frac{\alpha -\psi(\beta)}{\varphi(\alpha)-\beta}\frac{\beta}{\psi(\beta)},\end{equation}
and (iii)~that $\bar Y(T_\beta)- Y(T_\beta)$ is distributed as $-\underline Y(T_\beta)$ and that it is exponentially distributed with parameter $\psi(\beta)$. If the net cumulative claim process has a negative drift (i.e., ${\mathbb E}\, Y(1)<0$, in the ruin literature referred to as the {\it net profit condition}), then we can send $\beta$ to~$0$, so as to obtain
\[  \xi(\alpha)\coloneqq{\mathbb E}\,e^{-\alpha \bar Y(\infty)}=\frac{\alpha\,\varphi'(0)}{\varphi(\alpha)}.\]

In this paper we have two target objects. In the first place we wish to analyze $p(u,T_\beta\,|\,\Omega)$ through its transform with respect to the exceedance level $u$. Concretely, we point out, for various types of inter-inspection times $\Omega$, how one can compute
\[\pi(\alpha,\beta\,|\,\Omega)\coloneqq \int_0^\infty e^{-\alpha u} p(u,T_\beta\,|\,\Omega)\,{\rm d}u.\]
In the second place, we identify the Cram\'er-Lundberg asymptotics (in the regime that $u$ is large, that is) of the all-time bankruptcy probability $p(u\,|\,\Omega)\coloneqq p(u,\infty\,|\,\Omega)$, for which we also devise an efficient importance sampling based simulation approach.

\section{Exponential and hyperexponential inter-inspection times}\label{Exp_HypExp}
In this section we wish to study the bankruptcy probability ${p}(u,T_{\beta}\,|\,\Omega)$ for exponentially and hyperexponentially distributed inter-inspection times $\Omega$, by considering the corresponding transform $\pi(\alpha,\beta\,|\,\Omega)$. 
In the first subsection we consider the case of exponential inter-inspection times, i.e., we associate with $\Omega$ a random variable that is distributed as $T_\omega$ for some $\omega>0.$ As mentioned, this case has been dealt with already in \cite{BM2}, albeit relying on another argumentation, but we include it to maximally transparently demonstrate the approach followed in this paper, thus making the derivations for the other inter-inspection time distributions more accessible. 
In the second subsection we show how to use these results to deal with the case of hyperexponential inter-inspection times, in which we associate with $\Omega$ a random variable that is with probability $p_i>0$ distributed as $T_{\omega_i}$ for some $\omega_i>0$, where $\sum_{i=1}^d p_i=1$ for $d\in{\mathbb N}$.

In both cases, and also in the cases dealt with in the next section, our derivation consists of two steps. In the first step, we derive an equation for the transform $\pi(\alpha,\beta\,|\,\Omega)$, by conditioning on the evolution of the process until the first inspection time (if this occurs before killing). This equation, which is of a linear form, can be readily solved, but still contains one or multiple unknowns. In the second step, it is pointed out how these unknowns can be determined. 

\subsection{Exponentially distributed inter-inspection times}\label{Exp inter-inspection}
In this subsection we consider the case of inter-inspection times that are exponentially distributed with parameter $\omega>0$. 
Note that the time till the first event, which can be either an inspection or killing, is $T_{\beta+\omega}$ (i.e., exponentially distributed with parameter $\beta+\omega$).
We define the two non-negative random quantities 
\[Z^+ (\beta,\omega) \equiv Z^+ \coloneqq \bar{Y}(T_{\beta+\omega}), \quad Z^- (\beta,\omega) \equiv Z^- \coloneqq -\underline{Y}(T_{\beta+\omega}). 
\]
The random variable $Z^-$ is exponentially distributed with parameter $\theta\coloneqq \psi(\beta+\omega)$, and by virtue of the Wiener-Hopf decomposition $Z^-$ and $Z^+$ are independent. Below we use the compact notation $Z\coloneqq Z^+-Z^-$.
Note that, due to \eqref{xidef},
\begin{equation}
    \label{eq2.5}
\e \,e^{-\alpha Z^{+}} =\xi(\alpha,\beta+\omega)= \frac{\alpha -\psi(\beta+\omega)}{\varphi(\alpha)-\beta-\omega}\frac{\beta+\omega}{\psi(\beta+\omega)}.\end{equation}

\subsubsection*{-- Equation for transform}
We shall evaluate the quantity $\pi(\alpha,\beta\,|\,T_\omega)$ by conditioning on the first event. If $Y(\cdot)$ has to exceed level $u$ at an inspection time, then this first event time should be an inspection (i.e., not killing), which happens with probability $\omega/(\beta+\omega)$. Then there are two options: at this inspection time we are either already above $u$ and hence ruin occurs, or we are at some level $v$ below $u$ and we still need to bridge the remaining $u-v$. By following this line of reasoning, we obtain an expression in terms of our target quantity ${\pi}(\alpha,\beta\,|\,T_\omega)$, which we can solve. To this end, we first we write
\[{p}(u,T_\beta\,|\,T_\omega)=\frac{\omega}{\beta+\omega}\Big(\p(Z > u)+\int_{-\infty}^u \p(Z\in\diff v)\,{p}(u-v,T_\beta\,|\,T_\omega)\Big)
\]
In order to evaluate the associated transform, we multiply the entire equation by $e^{-\alpha u}$ and integrate over positive $u$:
\begin{align}
{\pi}(\alpha,\beta\,|\,T_\omega)=\frac{\omega}{\beta+\omega}\Bigg(&\underbrace{\int_0^\infty e^{-\alpha u}\,\p(Z >u)\diff u}_{\text{(i)}}+\nonumber\\
&\underbrace{\int_{u=0}^\infty e^{-\alpha u}\int_{v=-\infty}^u \p(Z\in\diff v)\,{p}(u-v,T_\beta\,|\,T_\omega)\diff u}_{\text{(ii)}}\Bigg)\label{eq1}
\end{align}
We proceed by separately analyzing the integrals (i) and (ii) appearing in \eqref{eq1}. 
\begin{itemize}
    \item[(i)] As a first step, we condition on the value of $Z^-$ and change the order of integration. With $\theta\coloneqq \psi(\beta+\omega)$, recalling that $Z^-$ is exponentially distributed with parameter $\theta$, we readily find
    \begin{align*}
        \int_0^\infty e^{-\alpha u}\ \p(Z > u)\diff u
        &=\theta\int_{u=0}^\infty e^{-\alpha u}\int_{v=0}^\infty e^{-\theta v}\p(Z^+> u+v)\diff v\diff u
    \end{align*}
    Performing the change of variable $x\coloneq u+v$, we have that this equals
    \begin{align*}
      \theta\int_{x=0}^\infty e^{-\theta x}\p(Z^+> x)\Big(\int_{u=0}^x e^{u(\theta-\alpha)}\diff u\Big)\diff x
    \end{align*}
    Upon evaluating the inner integral, this reduces to
    \begin{align*}
        \frac{\theta}{\theta-\alpha}\Big(\int_{0}^\infty e^{-\alpha x}\p(Z^+> x)\diff x -\int_{0}^\infty e^{-\theta x}\p(Z^+> x)\diff x\Big)
    \end{align*}
    Recalling that $\xi(\alpha,\beta+\omega)= \e \,e^{-\alpha Z^{+}}$, observe that
    \begin{align*}
        \int_0^\infty e^{-\alpha x} \p(Z^+ > x)\diff x
        &= \int_{x=0}^\infty e^{-\alpha x} \int_{y=x}^\infty \p(Z^+\in\diff y)\diff x\\
        &= \int_{y=0}^\infty \Big(\int_{x=0}^y e^{-\alpha x} \diff x\Big) \p(Z^+ \in \diff y)\\
        &= \frac{1}{\alpha}\Big(\int_{0}^\infty \p(Z^+\in\diff y) -\int_{0}^\infty e^{-\alpha y}\p(Z^+\in\diff y)\Big)\\
        &= \frac{1}{\alpha}\big(1-\xi(\alpha,\beta+\omega)\big)
    \end{align*}
    A similar reasoning applies when replacing $\alpha$ by $\theta$. Combining the above, we thus conclude that
    \begin{align*}
        \int_0^\infty e^{-\alpha u}\, \p(Z> u)\diff u 
        &=  \frac{\theta}{\theta-\alpha} \Big(\frac{1}{\alpha}\big(1-\xi(\alpha,\beta+\omega)\big)-\frac{1}{\theta}\big(1-\xi(\theta,\beta+\omega)\big)\Big)\\
        &=\frac{\alpha\xi(\theta,\beta+\omega)-\theta\xi(\alpha,\beta+\omega)}{\alpha(\theta-\alpha)}+\frac{1}{\alpha}\eqqcolon  F(\alpha,\beta, \omega).
    \end{align*}
    \item[(ii)] We continue by evaluating the second integral, starting by splitting it into two parts:
    \begin{align}
        \int_{u=0}^\infty e^{-\alpha u}\int_{v=-\infty}^u&\p(Z\in\diff v){p}(u-v,T_\beta\,|\,T_\omega)\,\diff u
        ={\pi}_1(\alpha,\beta\,|\,T_\omega)+{\pi}_2(\alpha,\beta\,|\,T_\omega),\label{eq2}
    \end{align}
    where the first term ${\pi}_1(\alpha,\beta\,|\,T_\omega)$ corresponds to the triangular integration area $\{0\leqslant v\leqslant u\}$ and the second term ${\pi}_2(\alpha,\beta\,|\,T_\omega)$ to the rectangular integration area $\{v\leqslant 0\leqslant u\}$. 
    
    \noindent As before, to evaluate the first double integral in \eqref{eq2}, we interchange the order of the integrals. We obtain that 
    \begin{align*}
        {\pi}_1(\alpha,\beta\,|\,T_\omega)={\pi}(\alpha,\beta\,|\,T_\omega)\int_0^\infty e^{-\alpha v}\p(Z\in\diff v),
    \end{align*}
    which, by conditioning on the value of $Z^-$, equals
    \begin{align*}
        \theta\ {\pi}(\alpha,\beta\,|\,T_\omega)\int_{v=0}^\infty e^{-\alpha v}\Big(\int_{x=0}^\infty e^{-\theta x}\p\big(Z^+-v\in\diff x\big)\Big)\diff v .
    \end{align*}
    By evaluating the inner integral with a change of variable $y\coloneq v+x$, this expression can be further rewritten as
    \begin{align*}
        \theta\ {\pi}(\alpha,\beta\,|\,T_\omega)&\int_{v=0}^\infty \int_{y=v}^\infty e^{v(\theta-\alpha)}e^{-\theta y}\p(Z^+\in\diff y)\diff v =\\
        &\theta\ {\pi}(\alpha,\beta\,|\,T_\omega)\int_{y=0}^\infty\Big(\int_{v=0}^y e^{v(\theta-\alpha)}\diff v\Big)e^{-\theta y}\p(Z^+\in\diff y),
    \end{align*}
   where we again changed the order of the integrals. Upon evaluating the obtained inner integral, this reduces to
    \begin{align*}
        \theta\ {\pi}(\alpha,\beta\,|\,T_\omega)\int_{0}^\infty e^{-\theta y} \Big(\frac{e^{\theta y-\alpha y}}{\theta-\alpha}-\frac{1}{\theta-\alpha}\Big)\,\p(Z^+\in\diff y).
    \end{align*}
    Some elementary manipulations lead us to the conclusion that
    \begin{align*}
        {\pi}_1(\alpha,\beta\,|\,T_\omega)&
        =\theta\ {\pi}(\alpha,\beta\,|\,T_\omega)\frac{\xi(\alpha,\beta+\omega)}{\theta-\alpha}-H(\alpha,\beta,\omega),
         \end{align*}
         where
         \[H(\alpha,\beta,\omega) \coloneqq \theta\ {\pi}(\alpha,\beta\,|\,T_\omega)\frac{\xi(\theta,\beta+\omega)}{\theta-\alpha}. \]

    \noindent In order to evaluate \eqref{eq2}, we are left with considering the second integral, i.e., ${\pi}_2(\alpha,\beta\,|\,T_\omega)$. Performing the change of variable $x\coloneq u-v$, and swapping the order of the integrals, we have that 
    \begin{align*}
        {\pi}_2(\alpha,\beta\,|\,T_\omega)&=\int_{v=-\infty}^0 \int_{x=-v}^\infty e^{-\alpha v}e^{-\alpha x}\ p(x,T_{\beta}\,|\,T_\omega)\,\diff x\,\p(Z\in\diff v)\\&=\int_{x=0}^\infty \int_{v=-x}^0 e^{-\alpha v}\p(Z\in\diff v)\,e^{-\alpha x}\ p(x,T_{\beta}\,|\,T_\omega)\diff x .
    \end{align*}
    Now note that the inner integral $v$ is only evaluated in non-positive values, so $-\infty<v\leqslant 0$. Since both $Z^-,Z^+\geqslant 0$, we have that $Z^+ -Z^-\leqslant 0$ if and only if $Z^-\geqslant Z^+$. The probability of $Z^+-Z^-$ being less than $v$ can hence be evaluated as follows: for $v\leqslant 0$, realizing that $\p(Z^-\geqslant Z^+)=\xi(\theta,\beta+\omega)$,
    \begin{align*}
        \p(Z<v)
        &=\p(Z<v\ |\ Z^-\geqslant Z^+)\,\p(Z^-\geqslant Z^+)=e^{\theta v}\ \xi(\theta,\beta+\omega),
    \end{align*}
    using the memoryless property of the exponential distribution. We thus find that
    \begin{align*}
        {\pi}_2(\alpha,\beta\,|\,T_\omega)=\theta\ \xi(\theta,\beta+\omega)\int_{x=0}^\infty\Big(\int_{v=-x}^0 e^{v(\theta-\alpha)}\diff v\Big) \,e^{-\alpha x} p(x,T_{\beta}\,|\, T_\omega)\diff x.
    \end{align*}
    Evaluating the (standard) inner integral and performing similar manipulations as before, we conclude that
    \begin{align*}
        {\pi}_2(\alpha,\beta\,|\,T_\omega)&=H(\alpha,\beta,\omega)-G(\alpha,\beta,\omega),
    \end{align*}
    where
    \[G(\alpha,\beta,\omega)\coloneqq \theta\ \xi(\theta,\beta+\omega)\frac{{\pi}(\theta,\beta\,|\,T_\omega)}{\theta-\alpha}\]
    Upon combining the above, we obtain that
    \[{\pi}_1(\alpha,\beta\,|\,T_\omega)+{\pi}_2(\alpha,\beta\,|\,T_\omega)=\theta\ {\pi}(\alpha,\beta\,|\,T_\omega)\frac{\xi(\alpha,\beta+\omega)}{\theta-\alpha}-G(\alpha,\beta,\omega).\]
\end{itemize}
We now use the results obtained above in order to derive an expression for the transform $ {\pi}(\alpha,\beta\,|\,T_\omega)$. Upon adding up the expressions that we found in (i) and (ii), we obtain
\begin{align}
  {\pi}(\alpha,\beta\,|\,T_\omega)=\frac{\omega}{\beta+\omega}\Bigg( F(\alpha,\beta,\omega)+\frac{\theta\,\xi(\alpha,\beta+\omega)}{\theta-\alpha}{\pi}(\alpha,\beta\,|\,T_\omega)-G(\alpha,\beta,\omega)\Bigg)\label{eq3}
\end{align}
Substituting the known expression for $\xi(\alpha,\beta+\omega)$, as was given in \eqref{eq2.5}, into  the relation \eqref{eq3}, we obtain the following result characterizing  ${\pi}(\alpha,\beta\,|\,T_\omega)$. 
\begin{proposition}\label{Expproposition} For any $\alpha\geqslant 0$ and $\beta>0$,
\begin{align}\label{explr}
    {\pi}(\alpha,\beta\,|\,T_\omega)\Big(1-\frac{\omega}{\beta+\omega-\varphi(\alpha)}\Big)=\frac{\omega}{\beta+\omega}\Big(F(\alpha,\beta,\omega)-G(\alpha,\beta,\omega)\Big)
\end{align}
\end{proposition}

\subsubsection*{-- Determination of the unknown constant} 
Whereas $F(\alpha,\beta,\omega)$ is an explicit expression in terms of the model primitives, $G(\alpha,\beta,\omega)$ contains the unknown quantity ${\pi}(\theta,\beta\,|\,T_\omega)$. This unknown can be identified by observing that any root  (in the right-half of the complex $\alpha$-plane, that is) of the left hand side of \eqref{explr} is necessarily also a root of the right hand side, using that ${\pi}(\alpha,\beta\,|\,T_\omega)$ is finite for all $\alpha,\beta>0$. Solving
\begin{align*}
    1=\frac{\omega}{\beta+\omega-\varphi(\alpha)},
    \end{align*}
leads us to $\varphi(\alpha)=\beta$, or $\alpha=\psi(\beta)$ since $\psi(\cdot)$ is the right-inverse of the Laplace exponent $\varphi(\cdot)$. It requires an elementary calculation to verify that we thus end up with
\begin{align*}
    {\pi}(\theta,\beta\,|\,T_\omega)=\frac{\omega\ \psi(\beta)\,\xi(\theta,\beta+\omega)-\beta\ \theta+\beta\ \psi(\beta)}{\theta\ \omega\ \psi(\beta)\,\xi(\theta,\beta+\omega)}.
\end{align*}
Plugging this expression into ${\pi}(\alpha,\beta\,|\,T_\omega)$, that we found in Proposition \ref{Expproposition}, after considerable calculus one arrives at the following result.

\begin{theorem}
For any $\alpha\geqslant 0$ and $\beta>0$,
\begin{align*}
    {\pi}(\alpha,\beta\,|\,T_\omega)=\frac{1}{\alpha}\frac{\omega\ \varphi(\alpha)(\theta-\alpha)\,\psi(\beta) +\alpha\ \beta \big(\varphi(\alpha)-\beta-\omega\big)(\theta-\psi(\beta))}{(\varphi(\alpha)-\beta)(\beta+\omega)(\theta-\alpha)\psi(\beta)}.
\end{align*}
\end{theorem}
Define $\varrho(\alpha,\beta\,|\,\Omega)\coloneqq  {\mathbb E}\,e^{-\alpha \bar Y_\Omega(T_\beta)}$, where 
\[\bar Y_\Omega(t)\coloneqq  \sup\{Y(\bar\Omega_n): \bar\Omega_n\leqslant t\}.\]

We note that ${\pi}(\alpha,\beta\,|\,T_\omega)$ and $\varrho(\alpha,\beta\,|\,T_\omega)$ can be expressed in one another, using the well-known `translation formula' $\varrho(\alpha,\beta\,|\,T_\omega)=1-\alpha\,{\pi}(\alpha,\beta\,|\,T_\omega)$. We eventually find the following result, which is a true generalization of the time-dependent version of the Pollaczek-Khinchine formula:
\begin{align*}
    \varrho(\alpha,\beta\,|\,T_\omega)=\frac{\alpha-\psi(\beta)}{\varphi(\alpha)-\beta}\frac{\beta}{\psi(\beta)}\frac{\varphi(\alpha)-\beta-\omega}{\alpha-\psi(\beta+\omega)}\frac{\psi(\beta+\omega)}{\beta+\omega};
\end{align*}
cf.\ \cite[Prop.\ 1]{BM2}.
Indeed, as $\omega\to\infty$, which corresponds to the inspection process taking place at an increasingly high frequency, thus approaching permanent inspection of the process $Y(t)$, we find that the transform corresponds to the one appearing in the time-dependent Pollaczek-Khinchine formula
\[\lim_{\omega\to\infty}\varrho(\alpha,\beta\,|\,T_\omega)= \frac{\alpha-\psi(\beta)}{\varphi(\alpha)-\beta}\frac{\beta}{\psi(\beta)};\]
see e.g.\ \cite[Thm.\ 4.1]{DebM}.

\subsection{Hyperexponentially distributed inter-inspection times}
We now consider the bankruptcy model with hyperexponential-$d$ inter-inspection times, by mimicking the line of reasoning that has been used in the previous subsection for the case of exponentially distributed inter-inspection times.

With probability $p_i> 0$ the inter-inspection time is exponentially distributed with rate $\omega_i>0$, for $i\in \{1,\dots,d\}$ and $\sum_{i=1}^d p_i=1$. 
This means that,  with ${\bs\omega}=(\omega_1,\ldots,\omega_d)$ and ${\bs p}=(p_1,\ldots,p_d)$, a generic inter-inspection time $T_{{\bs \omega,\bs p}}$ is characterized by the density,
\[{\mathbb P}(T_{{\bs \omega},{\bs p}}\in{\rm d}t) = \sum_{i=1}^d p_i\,\omega_ie^{-\omega_i t}{\rm d}t.\]
Recalling the definition of the two independent, non-negative random variables   $Z^+ (\beta,\omega_i) \coloneqq \bar{Y}(T_{\beta+\omega_i})$ and $Z^- (\beta,\omega_i) \coloneqq -\underline{Y}(T_{\beta+\omega_i})$, we find that
\[\xi_i(\alpha,\beta)\coloneqq \e \,e^{-\alpha Z^+(\beta,\omega_i)}= \xi(\alpha,\beta+\omega_i)=\frac{\alpha -\psi(\beta+\omega_i)}{\varphi(\alpha)-\beta-\omega_i}\frac{\beta+\omega_i}{\psi(\beta+\omega_i)},\]
and that $Z^- (\beta,\omega_i)$ is exponentially distributed with parameter $\theta_i\coloneqq \psi(\beta+\omega_i)$. 

\subsubsection*{-- Equation for transform}
The main idea is that at time $0$ we sample which of the $d$ exponentially distributed inter-inspection time will be used. Given this is the $i$-th inter-inspection time (with parameter $\omega_i$, that is), we can reuse the argument developed for the exponential case, but with $Z^+(\beta,\omega)$ replaced by $Z^+(\beta,\omega_i)$, and with $Z^-(\beta,\omega)$ replaced by $Z^-(\beta,\omega_i)$.
Hence, as the counterpart of the relation for the transform ${\pi}(\alpha,\beta\,|\,T_{\omega})$ that we derived in the case of exponentially distributed inter-inspection times, we now obtain
\begin{align*}
    {\pi}(\alpha,\beta\,|\,T_{{\bs \omega},{\bs p}})=&\ \sum_{i=1}^d p_i\ \frac{\omega_i}{\beta+\omega_i}\Bigg(\frac{\alpha\ \xi_i(\theta_i,\beta)-\theta_i\ \xi_i(\alpha,\beta)}{\alpha(\theta_i-\alpha)}+\frac{1}{\alpha}\:+\\\
    &\theta_i{\pi}(\alpha,\beta\,|\,T_{{\bs \omega},{\bs p}})\frac{\xi_i(\alpha,\beta)-\xi_i(\theta_i,\beta)}{\theta_i-\alpha}\:+\\
    &\theta_i\ \xi_i(\theta_i,\beta)\frac{{\pi}(\alpha,\beta\,|\,T_{{\bs \omega},{\bs p}})-{\pi}(\theta_i,\beta\,|\,T_{{\bs \omega},{\bs p}})}{\theta_i-\alpha}\Bigg);
\end{align*}
cf.\ the decomposition featuring in \eqref{eq3}.
Following the same procedure as in the case of exponentially distributed inter-inspection times, this can be rewritten as
\[ {\pi}(\alpha,\beta\,|\,T_{{\bs \omega},{\bs p}})=\ \sum_{i=1}^d p_i\ \frac{\omega_i}{\beta+\omega_i}\Big(F(\alpha,\beta,\omega_i)+\frac{\theta_i\,\xi_i(\alpha,\beta)}{\theta_i-\alpha}{\pi}(\alpha,\beta\,|\,T_{{\bs \omega},{\bs p}})-G_i(\alpha,\beta,{\bs \omega})\Big).
\]
Here the function $F(\alpha,\beta,\omega)$ is as defined in the subsection on exponentially distributed inter-inspection times, so that, for $i=1,\ldots,d$,
\[F(\alpha,\beta, \omega_i)=\frac{\alpha\xi(\theta_i,\beta+\omega_i)-\theta_i\xi(\alpha,\beta+\omega_i)}{\alpha(\theta_i-\alpha)}+\frac{1}{\alpha};\] in addition, for $i=1,\ldots,d$,
\[G_i(\alpha,\beta,{\bs \omega})\coloneqq  \bar G(\alpha,\beta,\omega_i)\,z_i,\]
where
\[\bar G(\alpha,\beta,\omega_i)\coloneqq \frac{\theta_i\ \xi_i(\theta_i,\beta)}{\theta_i-\alpha},\:\:\:\:
z_i\equiv z_i(\beta,{\bs\omega}) \coloneqq  {\pi}(\theta_i,\beta\,|\,T_{{\bs \omega},{\bs p}}).\]
We proceed by isolating our target object ${\pi}(\alpha,\beta\,|\,T_{{\bs \omega},{\bs p}})$. By inserting the expression for $\xi_i(\alpha,\beta)$, we thus end up with the following characterization.
\begin{proposition}\label{Hypexpproposition} For any $\alpha\geqslant 0$ and $\beta>0$,
\begin{align*}
    {\pi}(\alpha,\beta\,|\,T_{{\bs \omega},{\bs p}})\Big(1-\sum_{i=1}^d p_i \frac{\omega_i}{\beta+\omega_i-\varphi(\alpha)}\Big)=\sum_{i=1}^d p_i\frac{\omega_i}{\beta+\omega_i}\Big(F(\alpha,\beta,\omega_i)-G_i(\alpha,\beta,{\bs \omega})\Big).
\end{align*}
\end{proposition}

\subsubsection*{-- Determination of the unknown constants}
Observe that, in the right hand side of the equation featuring in Proposition \ref{Hypexpproposition}, the object $F(\alpha,\beta,\omega_i)$ is explicitly expressed in terms of the model primitives, and so is $\bar G(\alpha,\beta,\omega_i)$.
The constants $z_1,\ldots,z_d$, however, are yet to be determined. In the remainder of this subsection, we point out how these constants can be identified. 

\begin{lemma} \label{rootshypexp} For any $\beta$ and ${\bs\omega}$, the equation 
\begin{equation}\label{eq_target}\sum_{i=1}^d p_i \frac{\omega_i}{\beta+\omega_i-\varphi(\alpha)}=1\end{equation}
has $d$ roots in the right half plane, say, $\alpha_1\equiv \alpha_1(\beta,{\bs\omega}),\ldots,\alpha_d\equiv \alpha_d(\beta,{\bs\omega}).$
\end{lemma}

{\it Proof.} The main idea behind the proof (as well as other proofs later in the paper) is the following.
Consider a $d$-dimensional Markov additive process (MAP), say $(W(\cdot),J(\cdot))\equiv(W(t),J(t))_{t\geqslant 0}$; we refer to e.g.\ \cite[Section XI.2]{ASM2} for a definition. With $q_{ij}$, for $i,j=1,\ldots,d$, denoting the transition rates corresponding to the modulating continuous-time Markov chain $J(\cdot)$ on $\{1,\ldots,d\}$, and with $\varphi_1(\cdot),\ldots,\varphi_d(\cdot)$ Laplace exponents of $d$ L\'evy processes,
define the $(i,j)$th entry of the $(d\times d)$-matrix $M(\alpha,{\bs \eta})$ by
\[m_{ij}(\alpha,{\bs \eta}) \coloneqq  (\varphi_i(\alpha)-\eta_i)\,1\{i=j\}+q_{ij},\]
for a given componentwise positive vector ${\bs\eta}$. From \cite[Theorem 1 \& Remark 2.1]{IBM} we know that, if the $d$ L\'evy processes featuring in $(W(\cdot),J(\cdot))$ are no subordinators, then the equation ${\rm det}\, M(\alpha,{\bs\eta})=0$ has, for any componentwise positive vector ${\bs\eta}$, exactly $d$ roots in the right-half of the complex plane. To prove the lemma, the objective is now to construct a MAP $(W(\cdot),J(\cdot))$ such that the corresponding equation ${\rm det}\, M(\alpha,{\bs \eta})=0$ is, for some choice of the vector ${\bs\eta}$,  equivalent to Equation \eqref{eq_target}. 

Consider the  MAP defined by a transition rate matrix $Q$ that is given by
\begin{equation*}
    Q=\boldsymbol{p}\boldsymbol{\omega}^{\top}-\diag(\boldsymbol{\omega}),
\end{equation*}
and Laplace exponents $\varphi_1(\alpha)=\ldots=\varphi_d(\alpha)=\varphi(\alpha)$. The entries of the componentwise positive vector ${\bs\eta}$ are defined to be
$\eta_i=\beta$, for $i=1,\ldots,d$.
We now observe that
\begin{equation*}
M(\alpha,{\bs \eta})=Q+ {\rm diag}({\bs\varphi(\alpha)-\bs\eta})=\boldsymbol{p}\boldsymbol{\omega}^\top+{\rm diag}(\varphi(\alpha){\bs 1}-{\bs\omega}-{\beta}{\bs 1}),
\end{equation*}
with ${\bs 1}$ a $d$-dimensional all-ones vector. 
As a consequence, $M(\alpha,\eta)$ is the sum of a diagonal matrix and a rank-one matrix.

The last step is to evaluate ${\rm det}\,M(\alpha,{\bs \eta})$, to show that the equation ${\rm det}\,M(\alpha,{\bs \eta})=0$ is equivalent to Equation \eqref{eq_target}. To this end, define $\Delta\coloneqq {\rm diag}(\varphi(\alpha){\bs 1}-{\bs\omega}-{\beta}{\bs 1})$, and let $\delta_i$ be the $(i,i)$-th entry of this diagonal matrix. Observe that $M(\alpha,{\bs \eta})= \Delta+\boldsymbol{p}\boldsymbol{\omega}^\top$. As a consequence, we can write
\begin{equation*}
    \text{det}\,M(\alpha,{\bs \eta})=\text{det}\,\Delta\cdot\text{det}\,(I_d +\boldsymbol{p}\boldsymbol{\omega}^\top \Delta^{-1})
\end{equation*}
For $A$ of dimension $m\times n$, and $B$ of dimension $n\times m$, we have ${\rm det}\,(I_m+AB)={\rm det}\,(I_n+BA)$. We thus conclude that
\begin{align*}
    {\rm det}\,M(\alpha,{\bs \eta})&={\rm det}\,\Delta\cdot{\rm det}\,(I_1 + \bs\omega^\top \Delta^{-1}\boldsymbol{p})=\left(\prod_{i=1}^d \delta_i\right)\Bigg(\sum_{i=1}^d p_i\frac{\omega_i}{\varphi(\alpha)-\beta-\omega_i}+1\Bigg),
\end{align*}
so that indeed ${\rm det}\,M(\alpha,{\bs \eta})=0$ is equivalent to Equation \eqref{eq_target}.
\hfill$\Box$

\vb

To identify the unknown quantities $z_i$, we use the property that any root of the left hand side of the equation featuring in Proposition \ref{Hypexpproposition} (in the right-half of the complex $\alpha$-plane) is necessarily a root of the right hand side as well, noting that ${\pi}(\alpha,\beta\,|\,T_{{\bs \omega},{\bs p}})$ is finite for all $\alpha\geqslant 0$ and $\beta>0$. 
Let $\alpha_1,\ldots,\alpha_d$ be the roots of the left hand side in the right-half of the complex $\alpha$-plane; there are precisely $d$ of those, by Lemma \ref{rootshypexp}. 
Consequently we are to solve a system of linear equations: we have, for any $j=1,\ldots,d$,  
\begin{equation}\label{constz}
  \sum_{i=1}^d p_i\ \frac{\omega_i}{\beta+\omega_i }\Big(F(\alpha_j,\beta,\omega_i)-\bar G(\alpha_j,\beta,\omega_i)\,z_i\Big)=0.
\end{equation}
It can be argued that  the roots $\alpha_j$ are distinct; see Remark \ref{R1}. Hence, the unknown vector ${\bs z}$ can be identified by solving the $j$~linear equations \eqref{constz}. 
We thus arrive at the following result. 

\begin{theorem}
For any $\alpha\geqslant 0$ and $\beta>0$,
\begin{align*}
   {\pi}(\alpha,\beta\,|\,T_{{\bs \omega},{\bs p}})=\frac{\displaystyle \sum_{i=1}^d p_i\frac{\omega_i}{\beta+\omega_i}\Big(F(\alpha,\beta,\omega_i)-\bar G(\alpha,\beta, \omega_i)\,z_i\Big)}{\displaystyle 1-\sum_{i=1}^d p_i \frac{\omega_i}{\beta+\omega_i-\varphi(\alpha)}},
\end{align*}
where the vector ${\bs z}$ follows from the linear equations \eqref{constz}. 
\end{theorem}

\begin{remark} \label{R1} {\em 
The $d$ roots are actually real, as follows from a reasoning similar to the one in \cite[Section 3.2]{CDMR}. First realize that without loss of generality we can assume the rates $\omega_i$ are distinct (because if they were not, some of the $p_i$ could be `lumped' to make the $\omega_i$ distinct). Putting (without loss of generality) these $\omega_i$ in increasing order, we obtain that there is precisely one real root between $x_{i-1}$ and $x_i$, with $x_0\coloneqq 0$ and $x_i\coloneqq \psi(\beta+\omega_{i})$ for $i=1,\ldots,d$. Observe that the smallest root is $\psi(\beta)$.  This implies that the roots are distinct. \hfill$\Diamond$
}\end{remark}

\section{Erlang and hyper-Erlang inter-inspection times}\label{Erl_HypErl}

In this section we wish to study the bankruptcy probability $p(u,T_\beta\,|\,\Omega)$ for Erlang and hyper-Erlang distributed inter-inspection times $\Omega$, again by considering the corresponding transform $\pi(\alpha,\beta\,|\,\Omega)$. Interestingly, this will eventually enable us to approximate deterministically distributed inter-inspection times: observe that the sum of $k$ i.i.d.\ exponentially distributed random variables with parameter $\lambda k$ has mean $1/\lambda$ and variance $1/(k\lambda^2)$, with the latter decreasing to $0$ as $k\to\infty$.

In the first subsection we consider the Erlang distributed inter-inspection time case, where $\Omega$ is associated with a random variable $E_{k,\omega}$ that is the sum of $k\in{\mathbb N}$ independent random variables, each of them being distributed as $T_{\omega}$, for some $\omega>0$. In the second subsection we show how to use these results to handle the hyper-Erlang inter-inspection time case, in which we associate with $\Omega$ a random variable that is with probability $p_i$ distributed as $E_{k_{i},\omega}$ for some $k_i,\omega>0$, where $\sum_{i=1}^d p_i=1$. In both cases our derivation consists of the same two steps as in the previous section.

\subsection{Erlang distributed inter-inspection times}\label{EIIT}
We consider the bankruptcy model with Erlang-$k$ inspection times, with shape parameter $k\in{\mathbb N}$ and scale parameter $\omega$. This means that a generic inter-inspection time $E_{k,\omega}$ is characterized via the density
\[\p\big(E_{k,\omega}\in\diff t\big)=\omega\frac{(\omega t)^{k-1}}{(k-1)!}e^{-\omega t}\,{\rm d}t.
\]
In this subsection again a crucial role is played by the two non-negative random quantities
\[Z^+(\beta,\omega)\equiv Z^+\coloneqq \bar{Y}(T_{\beta+\omega}),\:\:\:Z^-(\beta,\omega)\equiv Z^-\coloneqq -\underline{Y}(T_{\beta+\omega}).\]
In what follows, let $Z^+_k\equiv Z^+_k(\beta,\omega)$ denote the sum of $k$ i.i.d.\ copies of $Z^+$, and analogously $Z^-_k\equiv Z^-_k(\beta,\omega)$  the sum of $k$ i.i.d.\ copies of $Z^-$.

Recalling ${\mathbb E}\,e^{-\alpha Z^+}=\xi(\alpha,\beta+\omega)$,  which can be evaluated by applying \eqref{xidef}, we write
\[
 \xi_k(\alpha,\beta+\omega)\coloneqq {\mathbb E}\,e^{-\alpha Z_k^+} =(\xi(\alpha,\beta+\omega))^k.\] In addition, we recall that $Z^-$ is exponentially distributed with parameter $\theta\coloneqq \psi(\beta+\omega)$. It is also noted that $Z^-$ and $Z^+$ are, by virtue of the Wiener-Hopf decomposition, independent. This entails that $Z_k^-$ is Erlang distributed with shape parameter $k$ and scale parameter $\theta$, and that $Z^-_k$ and $Z^+_k$ are independent. Below we use the compact notation $Z_k\coloneqq Z_k^+-Z_k^-$.

\subsubsection*{-- Equation for transform}

\noindent
The objective is to evaluate the transform
\[\pi(\alpha,\beta\,|\,E_{k,\omega})\coloneqq \int_0^\infty e^{-\alpha u} \, p(u,T_\beta\,|\,E_{k,\omega})\,{\rm d}u.\]
To this end, first notice that
\begin{equation*}
p(u,T_\beta\,|\,E_{k,\omega})=\left(\frac{\omega}{\beta+\omega}\right)^k\Big(\,{\mathbb P}(Z_k> u)+\int_{-\infty}^u
{\mathbb P}(Z_k\in {\rm d}v)\,p(u-v,T_\beta\,|\,E_{k,\omega})\,\Big).
\end{equation*}
In order to evaluate the associated transform, we multiply the entire equation by $e^{-\alpha u}$ and integrate over $u$, so as to obtain
\begin{align}\label{El1}
{\pi}(\alpha,\beta\,|\,E_{k,\omega})=\Big(&\frac{\omega}{\beta+\omega}\Big)^k\Bigg(\underbrace{\int_0^\infty e^{-\alpha u}\,\p(Z_k >u)\diff u}_{\text{(i)}}\:+\nonumber\\
&\underbrace{\int_{u=0}^\infty e^{-\alpha u}\int_{v=-\infty}^u \p(Z_k\in\diff v)\,{p}(u-v,T_\beta\,|\,E_{k,\omega})\diff u}_{\text{(ii)}}\Bigg)
\end{align}
As before, we deal with the terms (i) and (ii) separately. Lemma \ref{LL2} below, proven in the appendix, provides an expression for (i). 
Define, with $\xi^{(m)}(\alpha,\beta)$ denoting the $m$-th derivative of $\xi(\alpha,\beta)$ with respect to $\alpha$,
\begin{align}
I(\alpha,\beta,\omega)\coloneqq&\:\left(\frac{\theta}{\theta-\alpha}\right)^k \frac{1-\xi_k(\alpha,\beta+\omega)}{\alpha}-\frac{1}{\theta} \sum_{n=0}^{k-1}\left(\frac{\theta}{\theta-\alpha}\right)^{k-n}
\left(1-\sum_{m=0}^n \xi_k^{(m)}(\theta,\beta+\omega)\frac{(-\theta)^m}{m!}\right).\label{ErlangInt1}
\end{align}
\begin{lemma} \label{LL2} For any $\alpha\geqslant 0$,
\begin{align*}
\int_0^\infty e^{-\alpha u}\,{\mathbb P}(Z_k> u)\,{\rm d}u=&\:I(\alpha,\beta,\omega).\label{ErlangInt1}
\end{align*}
\end{lemma}
The right hand side of \eqref{ErlangInt1} can be further simplified by swapping the order of summation and evaluating the finite geometric series, but we leave this computation out.  

We now consider integral (ii). To this end, define 
\[\delta_{n,k}\coloneqq \frac{(-\theta)^n}{n!}\xi_k^{(n)}(\theta,\beta+\omega).
\]
Also,
\begin{align*}
 J(\alpha,\beta,\omega)&\coloneqq  \sum_{n=0}^{k-1}\delta_{k-1-n,k}\left(\frac{\theta}{\theta-\alpha}\right)^{n+1}
\sum_{m=0}^n \frac{(\alpha-\theta)^m}{m!}\pi^{(m)}(\theta,\beta\,|\,E_{k,\omega})=
\sum_{i=0}^{k-1} \bar J_i(\alpha,\beta,\omega)\,z_i,
\end{align*}
where (by swapping the order of the two sums) $J_i(\alpha,\beta,\omega)$ and $z_i\equiv z_i(\beta,\omega)$ are
\begin{equation} \label{defjbar}\bar J_i(\alpha,\beta,\omega)\coloneqq \sum_{n=i}^{k-1} \delta_{k-1-n,k}\left(\frac{\theta}{\theta-\alpha}\right)^{n+1}\frac{(\alpha-\theta)^i}{i!},\:\:\:\:
z_i\coloneqq \pi^{(i)}(\theta,\beta\,|\,E_{k,\omega}).\end{equation}
Observe that $\bar J_i(\alpha,\beta,\omega)$ is explicitly known (in terms of the model primitives, that is), whereas the $k$ coefficients $z_0,\ldots,z_{k-1}$ are yet to be determined.  The following lemma is proven in the appendix.

\begin{lemma} \label{LL3} For any $\alpha\geqslant 0$,
\begin{align}
\int_{u=0}^\infty e^{-\alpha u}\int_{v=-\infty}^u &\p(Z_k\in\diff v)\,{p}(u-v,T_\beta\,|\,E_{k,\omega})\diff u\nonumber\\
&=-J(\alpha,\beta,\omega)+\Big(\frac{\theta}{\theta-\alpha}\Big)^k\xi_k(\alpha,\beta+\omega)\,\pi(\alpha,\beta\,|\,E_{k,\omega}).\nonumber
\end{align}
\end{lemma}

We now use the results obtained above  to derive an expression for the transform $\pi(\alpha,\beta\,|\,E_{k,\omega})$. Adding up the expressions that we have found in Lemmas \ref{LL2} and \ref{LL3},
\begin{align}\label{Erlangdecomposition}
\pi(\alpha,\beta\,|\,E_{k,\omega})=\Big(\frac{\omega}{\beta+\omega}\Big)^k\bigg(&I(\alpha,\beta,\omega)-J(\alpha,\beta,\omega)\:+\nonumber\\
&\Big(\frac{\theta}{\theta-\alpha}\Big)^k\xi_k(\alpha,\beta+\omega)\,\pi(\alpha,\beta\,|\,E_{k,\omega})\bigg).
\end{align}
Noticing that 
\[\left(\frac{\theta}{\theta-\alpha}\right)^k
\xi_k(\alpha,\beta+\omega)= \left(\frac{\beta+\omega}{\beta+\omega-\varphi(\alpha)}\right)^k, \]
we arrive at the following result.
\begin{proposition}\label{Erlangproposition}
For any $\alpha\geqslant 0$ and $\beta>0$,
\[\pi(\alpha,\beta\,|\,E_{k,\omega})\left(1-\left(\frac{\omega}{\beta+\omega-\varphi(\alpha)}\right)^k\right) =\Big(\frac{\omega}{\beta+\omega}\Big)^k \big(I(\alpha,\beta,\omega)-J(\alpha,\beta,\omega)\big).\]
\end{proposition}

\subsubsection*{-- Determination of the unknown constants}
We observe that in Proposition \ref{Erlangproposition} the right hand side of the equation contains the object $J(\alpha,\beta,\omega)$, which contains the unknown quantities $z_0,\ldots,z_{k-1}$. The aim of the remainder of this subsection is to identify these unknown quantities, leading to an expression for our object of study, i.e., ${\pi}(\alpha,\beta\,|\,E_{k,\omega})$.

\begin{lemma}\label{rootserlang}
For any $\beta,\,\omega$ and $k$, the equation
\begin{equation}\label{eqlemma3}
\left(\frac{\omega}{\beta+\omega-\varphi(\alpha)}\right)^k=1
\end{equation}
has $k$ roots in the right half of the complex plane,  say, $\alpha_1\equiv \alpha_1(\beta,{\omega}),\ldots,\alpha_k\equiv \alpha_k(\beta,{\omega}).$
\end{lemma}

{\it Proof.} We use the precise same line of reasoning as in the proof of Lemma \ref{rootshypexp}.
Again  $\varphi_1(\alpha)=\ldots=\varphi_k(\alpha)=\varphi(\alpha)$ and
$\eta_i=\beta$ for $i=1,\ldots,k$, but now we consider the following MAP, with the entries of the transition rate matrix $Q$ given by the cyclic structure
\[q_{i,i}=-\omega,\quad q_{i,i+1}=\omega,\quad q_{k,1}=\omega,
\]
for $i=1,\dots, k$. Hence,
\[M(\alpha,\bs\eta)= Q + \text{diag}\,(\bs {\varphi}(\alpha)-\bs {\eta}),
\]
which, unfortunately, does not have the nice property of being the sum of a diagonal matrix and a rank-one matrix (as was the case in the setting of the proof of Lemma \ref{rootshypexp}), since the matrix $Q$ is not of rank one. To show that the equation $\text{det}\, M(\alpha,\bs\eta)=0$ is equivalent to Equation \eqref{eqlemma3}, we express the determinant in terms of determinants of smaller matrices, by performing a Laplace expansion along the first row. To this end, define $\omega^\star\coloneqq \beta+\omega-\varphi(\alpha)$ and the matrices $\hat Q$ and $\check Q$ given by the entries $\hat q_{ii}=-\omega^\star,\,\hat q_{i,i+1}=\omega$ and $\check q_{ii}=\omega,\,\check q_{i,i-1}=-\omega^\star$, respectively. As a consequence, we can write 
\begin{align*}
    \text{det }M(\alpha,\bs\eta)&=-\omega^\star \text{ det}\, \hat Q + (-1)^{k+1}\omega\text{ det}\,\check Q\\
    &=(-\omega^\star)^k-(-\omega)^k
    =(-\omega+\varphi(\alpha)-\beta)^k-(-\omega)^k,
\end{align*}
so that $\text{det }M(\alpha,\bs\eta)=0$ is equivalent to Equation \eqref{eqlemma3}. This settles the claim. 
\hfill$\Box$

\vb

The unknowns $z_i$ can be found relying on Lemma \ref{rootserlang}, once more using that any root of the left hand side  (in the right-half of the complex $\alpha$-plane) of the equation featuring in Proposition \ref{Erlangproposition} is also a root of the right hand side. It entails that we have to solve the linear system
\begin{equation}
    \label{constz2}
I(\alpha_j,\beta,\omega) = \sum_{i=0}^{k-1} \bar J_i(\alpha_j,\beta,\omega)\,z_i,\end{equation}
for $j=1,\ldots,k.$ Below we write $J(\alpha,\beta,\omega):=\sum_{i=0}^{k-1} \bar J_i(\alpha,\beta,\omega)\,z_i.$

\begin{theorem}
For any $\alpha\geqslant 0$ and $\beta>0$,
\begin{align*}
   {\pi}(\alpha,\beta\,|\,E_{k,\omega})=
   \frac{\displaystyle  \left(\frac{\omega}{\beta+\omega}\right)^k\big(I(\alpha,\beta,\omega)-J(\alpha,\beta,\omega)\big)}{1-\left({\displaystyle \frac{\omega}{\beta+\omega-\varphi(\alpha)}}\right)^k},
\end{align*}
where the vector ${\bs z}$ follows from the linear equations \eqref{constz2}. 
\end{theorem}

\begin{remark} \label{R2} {\em The entries of ${\bs z}$ are real, as can be seen as follows. (i) First observe that the roots $\alpha_j$ solves the equation, with ${\rm i}$ denoting the imaginary unit,
\begin{equation}\label{complc}\varphi(\alpha) =\omega -\omega \,e^{2\pi {\rm i}j/k} +\beta,\end{equation}
with $j=0,\ldots,k-1.$ This means that any $\varphi(\alpha_j)$ is either real or the complex conjugate of $\varphi(\alpha_{j'})$ for $j\not=j'$. (ii) Noting that, for $v,w\in{\mathbb R}$, 
\[\varphi(v+w{\rm i})=\log {\mathbb E}e^{-(v+w{\rm i})Y(1)} = \log {\mathbb E} \left[e^{-vY(1)}(\cos(-wY(1))+{\rm i}\sin(-wY(1))\right],\]
it requires a straightforward verification to conclude that if $\alpha_j=v+w{\rm i}$ solves \eqref{complc}, then so does its complex conjugate $\alpha_{j'}=v-w{\rm i}$. Hence the roots of \eqref{eqlemma3} are either real, or belong to a pair of complex conjugates. (iii)~Consider row $j$ and $j'$ in the linear system \eqref{constz2}, with $\alpha_j$ and $\alpha_{j'}$ being complex conjugates. From their definitions it can be verified that $I(\alpha_j,\beta,\omega)$ is the complex conjugate of $I(\alpha_{j'},\beta,\omega)$, and that $\bar J_i(\alpha_j,\beta,\omega)$ is the complex conjugate of $\bar J_i(\alpha_{j'},\beta,\omega)$. Now replace rows $j$ and $j'$ by one row that follows by adding the two rows up, and a second row that follows by subtracting one row from the other, and dividing by ${\rm i}$. All coefficients in the new linear system have become real, and hence the solution ${\bs z}$ is real. \hfill$\Diamond$
}\end{remark}

\subsection{Hyper-Erlang distributed inter-inspection times}\label{HED}
We consider the bankruptcy model with hyper-Erlang-$d$ inter-inspection times, with common scale parameter. With probability $p_i$ the inter-inspection time is Erlang-$k_i$ distributed (for some $k_i\in{\mathbb N}$), with scale parameter $\omega>0$, for $i\in\{1,\dots,d\}$ and $\sum_{i=1}^d p_i=1$. This means that a generic inter-inspection time $E_{\,{\bs k},{\bs p},\omega}$ has the density, with ${\bs k}=(k_1,\dots, k_d)$ and ${\bs p}=(p_1,\dots,p_d)$,
\[\p(E_{\,{\bs k},{\bs p},\omega}\in \diff t)=\sum_{i=1}^d p_i \,\omega\frac{(\omega t)^{k_i -1}}{(k_i-1)!}e^{-\omega t}\,{\rm d}t.
\]
Without loss of generality, we assume that the $k_i$ are put in increasing order.

\begin{remark}{\em 
    In this subsection we work with a common scale parameter $\omega>0$. At the expense of more complex computations and the introduction of additional notation, we could have dealt with a setting in which with probability $p_i$ the inter-inspection time is Erlang-$k_i$ distributed with scale parameter $\omega_i>0$ (i.e., with the parameters $\omega_i$ not necessarily identical). We have chosen to work with a common scale parameter, as the proof of \cite[Theorem III.4.2]{ASM2} reveals that any distribution on the positive half-line can be approximated arbitrarily closely by mixtures of Erlang distributions with the same scale parameter. \hfill$\Diamond$
}\end{remark}

\subsubsection*{-- Equation for transform}
Based on the decomposition that we developed in \eqref{Erlangdecomposition}, corresponding to the transform $\pi(\alpha,\beta\,|\,E_{k,\omega})$ in the case of Erlang distributed inter-inspection times, we now obtain
\begin{align}\label{Hypexptransform}
\pi(\alpha,\beta\,|\,E_{\,{\bs k},{\bs p},\omega})=\sum_{i=1}^d p_i \Big(\frac{\omega}{\beta+\omega}\Big)^{k_i}\bigg(&I_i(\alpha,\beta,\omega)-J_i(\alpha,\beta,\omega)\:+\nonumber\\
&\Big(\frac{\theta}{\theta-\alpha}\Big)^{k_i}\xi_{k_i}(\alpha,\beta+\omega)\,\pi(\alpha,\beta\,|\,E_{\,{\bs k},{\bs p},\omega})\bigg), 
\end{align}
with
\begin{align*}
    I_i(\alpha,\beta,\omega)&\coloneqq\left(\frac{\theta}{\theta-\alpha}\right)^{k_i} \frac{1-\xi_{k_i}(\alpha,\beta+\omega)}{\alpha}\:-\\
&\hspace{2cm}\frac{1}{\theta} \sum_{n=0}^{{k_i}-1}\left(\frac{\theta}{\theta-\alpha}\right)^{{k_i}-n}
\left(1-\sum_{m=0}^n \xi_{k_i}^{(m)}(\theta,\beta+\omega)\frac{(-\theta)^m}{m!}\right),\\
J_i(\alpha,\beta,\omega)&\coloneqq\sum_{n=0}^{k_{i}-1}\delta_{{k_i}-1-n,k_i}\left(\frac{\theta}{\theta-\alpha}\right)^{n+1}
\sum_{m=0}^n \frac{(\alpha-\theta)^m}{m!}\pi^{(m)}(\theta,\beta\,|\,E_{\,{\bs k},{\bs p},\omega})=\sum_{j=0}^{k_{i}-1}\bar{J}_{i,j}(\alpha,\beta,\omega)\,z_j,
\end{align*}
where, again by swapping the order of the two sums
\begin{align*}
    \bar{J}_{i,j}(\alpha,\beta,\omega)\coloneqq \sum_{n=j}^{k_{i}-1}\delta_{k_{i}-1-n,k_i}\left(\frac{\theta}{\theta-\alpha}\right)^{n+1}\frac{(\alpha-\theta)^j}{j!},
\end{align*}
and $z_j\equiv z_j\,(\beta,\omega)\coloneqq \pi^{(j)}(\theta,\beta\,|\,E_{\,{\bs k},{\bs p},\omega}).$ Observe that $\bar J_{i,j}(\alpha,\beta,\omega)$ is explicitly known, whereas the coefficients $z_0,\dots,z_{k_{d}-1}$ are yet to be determined.
As before, the next step is to isolate the object $\pi(\alpha,\beta\,|\,E_{\,{\bs k},{\bs p},\omega})$ in Equation \eqref{Hypexptransform}. We thus end up with the following characterization of our target quantity.
\begin{proposition}\label{HyperErlangprop}
For any $\alpha\geqslant 0$ and $\beta>0$,
\begin{align*}
    \pi(\alpha,\beta\,|\,E_{\,{\bs k},{\bs p},\omega})&\Bigg(1-\sum_{i=1}^d p_i\Big(\frac{\omega}{\beta+\omega-\varphi(\alpha)}\Big)^{k_i}\Bigg)\\
    &=\sum_{i=1}^d p_i\Big(\frac{\omega}{\beta+\omega}\Big)^{k_i}\big(I_i(\alpha,\beta,\omega)-J_i(\alpha,\beta,\omega)\big).
\end{align*}
\end{proposition}

\subsubsection*{-- Determination of the unknown constants}
In the right-hand side of Proposition \ref{HyperErlangprop} we have that the object $I_i(\alpha,\beta,\omega)$ is explicitly expressed in terms of the model primitives, and so is $\bar J_i(\alpha,\beta,\omega)$. The coefficients $z_0,\dots,z_{k_{d}-1}$, however, are still unknown. The aim of the remainder of this subsection is to identify these unknown quantities. 

\begin{lemma}\label{eqlemma4++}
For any $\beta,\,\omega, {\bs p}$ and ${\bs k}$, the equation
\begin{equation}\label{eqlemma4}
\sum_{i=1}^d p_i\Big(\frac{\omega}{\beta+\omega-\varphi(\alpha)}\Big)^{k_i}=1
\end{equation}
has $k_d$ roots in the right half of the complex plane,  say, $\alpha_1\equiv \alpha_1(\beta,\bs{\omega}),\ldots,\alpha_{k_d}\equiv \alpha_{k_d}(\beta,\bs{\omega}).$
\end{lemma}

{\it Proof.} We use a similar reasoning as in the proofs of Lemmas \ref{rootshypexp} and \ref{rootserlang}. Without loss of generality we had assumed that the shape parameters ${\bs k}=(k_1,\dots, k_d)$ are ordered, such that $k_1<\dots<k_d$, but in addition we can assume $k_i=i$. Indeed, the density of a generic inter-inspection time $E_{\,{\bs k},{\bs p},\omega}$ can be written as 
\[\p(E_{\,{\bs k},{\bs p},\omega}\in\diff t)=\sum_{i=1}^{k_d} \bar p_i\, \omega \frac{(\omega t)^{i-1}}{(i-1)!}e^{-\omega t}\diff t,
\]
where $\bar p_{k_i} = p_{i}$ (for $i=1,\ldots,d$) and $0$ else. 
Hence, we wish to prove that
\begin{equation}\label{eqlemma4+}
\sum_{i=1}^{k_d} \bar p_i\Big(\frac{\omega}{\beta+\omega-\varphi(\alpha)}\Big)^{i}=1
\end{equation}
has precisely $k_d$ roots in the right half of the complex $\alpha$-plane.

Consider the MAP characterized by   $\varphi_1(\alpha)=\ldots=\varphi_{k_d}(\alpha)=\varphi(\alpha)$ and
$\eta_i=\beta$ for $i=1,\ldots,k_d$, while the entries of the first column of the transition rate matrix $Q$ given by, with $\mathfrak{p}_i\coloneqq \sum_{k=1}^{i}\bar p_i$ for $i=1,\ldots,k_d$,
\[q_{1,1}=-(1-\bar p_1)\,\omega,\quad q_{k_d,1}=\omega, \quad q_{i,1}=\frac{\bar p_i}{1-\mathfrak{p}_{i-1}}\omega \quad\text{ for }i=2,\dots,k_d-1,
\]
and the other entries given by
\[q_{i,i}=-\omega\quad \text{ for } i=2,\dots,k_d,\quad q_{i,i+1}=\frac{1-\mathfrak{p}_{i}}{1-\mathfrak{p}_{i-1}}\,\omega \quad \text{ for } i=1,\dots,k_d-1.
\]
As a result,
\[M(\alpha,\bs\eta)= Q + \text{diag}(\bs {\varphi}(\alpha)-\bs {\eta}).
\]
The rationale behind this choice of $Q$ is the following. Inspection of the proofs of Lemmas \ref{rootshypexp} and~\ref{rootserlang} reveals that $Q$ should `encode' the distribution of $\Omega$. In the present hyper-Erlang case, there are $i$ phases with probability $\bar p_i$, which means that with probability $(1-\mathfrak{p}_{i})/({1-\mathfrak{p}_{i-1}})$ one proceeds to the $(i+1)$-st phase after having finished the $i$-th phase (and with the complementary probability $\bar p_i/({1-\mathfrak{p}_{i-1}})$ one is done). 

To show that the equation $\text{det}\, M(\alpha,\bs\eta)=0$ is equivalent to Equation \eqref{eqlemma4}, we evaluate the determinant by Laplace expansion along the first column. To this end, we define $\omega^\star\coloneqq \beta+\omega-\varphi(\alpha)$. As a consequence, we can write
\begin{align*}
    \text{det}\,M(\alpha,\bs\eta)&=(\bar p_1\omega-\omega^\star)(-\omega^\star)^{k_d-1}+\sum_{j=2}^{k_d}(-1)^{j+1}\frac{\bar p_j}{1-\mathfrak{p}_{j-1}}\omega\prod_{i=1}^{j-1}\frac{1-\mathfrak{p}_{i}}{1-\mathfrak{p}_{i-1}}\omega^{j-1}(-\omega^\star)^{k_d-j}\\
    &=(\bar p_1\omega-\omega^\star)(-\omega^\star)^{k_d-1}-\sum_{j=2}^{k_d} \bar p_j\,(-\omega)^j(-\omega^\star)^{k_d-j}\\
    &=(-\omega^\star)^{k_d}-\sum_{j=1}^{k_d} \bar p_j\,(-\omega)^j(-\omega^\star)^{k_d-j},
\end{align*}
so that $\text{det}\,M(\alpha,\bs\eta)=0$ is equivalent to Equation \eqref{eqlemma4+}, as desired.
\hfill$\Box$

\vb

The unknowns $z_1,\ldots,z_{k_d}$ can be found using Lemma \ref{eqlemma4++} in combination with the property that any root (in the right half of the complex $\alpha$-plane) of the left hand side of the equation featuring in Proposition \ref{HyperErlangprop}, is necessarily a root of the right hand side as well. This reasoning leads to the linear system
\begin{equation}\label{constz3}
    \sum_{i=1}^d p_i \Big(\frac{\omega}{\beta+\omega}\Big)^{k_i}I_i(\alpha_\ell,\beta,\omega)=\sum_{i=1}^d p_i \Big(\frac{\omega}{\beta+\omega}\Big)^{k_i}\sum_{j=0}^{k_{i}-1}\bar{J}_{i,j}(\alpha_\ell,\beta,\omega)\,z_j.
\end{equation}
for $\ell=1,\dots,k_d$. As in Remark \ref{R2}, the entries of ${\bs z}$ are real. Defining
\[J_i(\alpha,\beta,\omega):=\sum_{j=0}^{k_{i}-1}\bar{J}_{i,j}(\alpha,\beta,\omega)\,z_j,\]
we have found the following result. 

\begin{theorem}
For any $\alpha\geqslant 0$ and $\beta>0$,
\begin{equation*}
    \pi(\alpha,\beta\,|\,E_{\,{\bs k},{\bs p},\omega})=\frac{\displaystyle \sum_{i=1}^d p_i \Big(\frac{\omega}{\beta+\omega}\Big)^{k_i}\Big(I_i(\alpha,\beta,\omega)-{J}_{i}(\alpha,\beta,\omega)\Big)}{\displaystyle 1-\sum_{i=1}^d p_i \Big(\frac{\omega}{\beta+\omega-\varphi(\alpha)}\Big)^{k_i}}
\end{equation*}
where the vector $\bs z$ follows from the linear equations \eqref{constz3}.
\end{theorem}

\section{Reflections on distribution of inter-inspection times}\label{refl}
In this section we discuss the role played by the inter-inspection times $\Omega.$
While the first subsection points out that phase-type distributions are a natural candidate to approximate $\Omega$ with, the second subsection argues that only a subclass of these phase-type distributions  (including the ones featuring in Sections \ref{Exp_HypExp} and \ref{Erl_HypErl}) allows an explicit analysis.

\subsection{Approximating general distributions, two-moments fit}\label{TMF}
The inter-inspection time distributions we have dealt with in the previous sections are all contained in the wider class of {\it phase-type distributions}. Informally, a phase-type random variable has a distribution that is constructed by taking convolutions and/or mixtures of exponential distributions. Phase-type distributions owe their popularity to the fact that they are capable of approximating any distribution (on the positive half-line, that is) arbitrarily closely \cite[Theorem III.4.2]{ASM2}, while still allowing (semi-)explicit evaluation. 

A phase-type random variable is uniquely characterized as the time until absorption of a finite-state Markov chain with a single absorbing state. Following the approach presented in  \cite[Section III.4]{ASM2}, it is encoded by an initial distribution ${\bs a}\in {\mathbb R}^m$ and a sub-stochastic $m\times m$ transition rate matrix $S$, where `sub-stochastic' means that in at least one row the sum of the non-diagonal elements exceeds the diagonal element. For instance, the Erlang distribution with $k$ phases is represented by $a_1=1$, $S_{ii}=-\omega$ for $i =1,\ldots,k$ and $S_{i,i+1}=\omega$ for $i=1,\ldots,k-1$ (while all other entries of the $k$-dimensional vector ${\bs a}$ and of the $k\times k$ matrix $S$ are set equal to $0$). 

It is noted that the proof of \cite[Theorem III.4.2]{ASM2}  reveals that actually any distribution on the positive half-line can be approximated arbitrarily closely by elements from a specific subclass of phase-type distributions. This smaller class are the mixtures of Erlang distributions, each of them having the same scale parameter, i.e., the class of distributions dealt with in Section \ref{HED}. The proof also indicates that the drawback of working with this specific class of phase-type distributions, is that one may need a large dimension $m$ to get an accurate fit. That is the reason why a two-moment fit has been proposed in which a relatively low value of $m$ suffices. In the one presented in \cite{TIJ}, for distributions with a coefficient of variation less than 1 a mixture of Erlang random variables (with he same scale parameter) is used, and for distributions with a coefficient of variation larger than 1 a hyperexponential random variable, all of them of relatively low dimension; this fitting method is compactly summarized in e.g.\ \cite[Section 3.1]{KUIP}. This two-moment fit renders the inter-inspection time distributions covered by Sections \ref{Exp_HypExp} and \ref{Erl_HypErl} particularly useful.

\subsection{Limitations}\label{limi}
In the previous sections we succeeded in identifying the transform  $\pi(\alpha,\beta\,|\,\Omega)$ for (hyper-)exponential and (hyper-)Erlang inter-inspection times, both belonging to the family of phase-type distributions. This  suggests that the analysis may extend to the full class of {phase-type} distributions. Careful inspection, however, reveals that such an extension is challenging.  

To explain what complications arise, let us analyze the distribution of $Z=Y(\Omega)$ if the inter-inspection times $\Omega$ are of phase-type. Let $K_j$, for $j=1,\ldots,m$, be the number of visits of the associated underlying Markov chain to state $j$ before absorption. Then $Z=Z^+-Z^-$, where, in self evident notation,
\[Z^+\coloneqq   \sum_{j=1}^m \sum_{i=1}^{K_j} Z_{ji}^+,\:\:\:Z^-\coloneqq \sum_{j=1}^m \sum_{i=1}^{K_j}Z_{ji}^{-}. \]
While all pair of random variables  $Z_{ji}^+$ and $Z_{ji}^{-}$ are independent, for obvious reasons
\[\sum_{i=1}^{K_j}Z_{ji}^+\:\:\:\mbox{and}\:\:\:\sum_{i=1}^{K_j}Z_{ji}^-\]
are {\it not}. As a consequence, for general phase-type $\Omega$ it is not clear how to evaluate objects like
\[  \int_0^\infty e^{-\alpha u}\ \p(Z > u)\diff u\:\:\:\mbox{and}\:\:\:{\mathbb P}(Z<v\,|\,Z<0),\]
which were the key quantities in the computations of the previous sections. 

As is readily verified, the reason why the computations in Sections \ref{Exp_HypExp} and \ref{Erl_HypErl} worked out is that, for all of the inter-inspection time distributions considered, the random variables $K_j$ are either 0 or~1, so that $Z^+$ and $Z^-$ are independent. A sufficient condition for this property to hold, is that the sub-stochastic transition rate matrix $S$ does not allow cycles (so that every state is visited at most once), a condition that is met for (hyper-)exponential and (hyper-)Erlang inter-inspection times.

\begin{remark}{\em 
We conclude this section by remarking that it {\it is} true that $Z^-$ is of phase-type, owing to the facts that (i) all random variables $Z_{ji}^-$ are exponentially distributed and that (ii) the random variables $K_j$, for $j=1,\ldots,m$, correspond to the number of visits of a Markov chain to its respective states (before absorption, that is). 

A second observation is that the expressions found in the previous sections suggest a relation with associated queueing systems; for instance, the expressions found for the model of Section~\ref{EIIT} strongly resemble those of the workload in the E$_{k}$/G/1 queue (i.e., the queue with Erlang-$k$ interarrival times, general service times and a single server); cf.\  \cite[Thm.\ 25, p.\ 44]{Prabhu}. One would therefore expect that analyzing our system with phase-type $\Omega$ is of the same complexity as the workload in a Ph/G/1 queue (i.e., the queue with phase-type interarrival times, general service times and a single server). This is, however, {\it not} true, as in the Ph/G/1 queue there is independence between the sequence of interarrival times and the sequence of service times, while in our bankruptcy context $Z^-$ and $Z^+$ are not necessarily independent. \hfill$\Diamond$}
\end{remark}

\section{Cram\'er-Lundberg-type asymptotics}\label{CLA}
The objective of this section is to find the exact asymptotics of 
\[p(u\,|\,\Omega)\coloneqq  {\mathbb P}(\exists n\in{\mathbb N}: Y(\bar\Omega_n)>u),\]
i.e., we aim to identify an explicit function $f(u)$ such that $f(u)\, p(u\,|\,\Omega)\to 1$ as $u\to\infty.$ We focus on the light-tailed regime, in which $p(u\,|\,\Omega)$ asymptotically looks like $\gamma e^{-\theta^\star u}$ for a positive decay rate $\theta^\star$ and `preconstant' $\gamma.$ This extends the classical Cram\'er-Lundberg ruin probability asymptotics to our bankruptcy setting.

In the first subsection we derive a useful identity, relating $p(u\,|\,\Omega)$ to the likelihood ratio of an exponentially-twisted version of the process $Y(\cdot)$ at a stopping time, revealing the value of $\theta^\star$. This relation is used in the second subsection: by adapting the techniques that we developed earlier in this paper, we succeed in identifying $\gamma$ for the types of inter-inspection time distributions considered in Sections~\ref{Exp_HypExp} and~\ref{Erl_HypErl}. 

\subsection{Change of measure}
We start by introducing some useful notation. 
As before, we consider the net cumulative claim process $Y(t)$ at the inspection epochs $\bar\Omega_1, \bar\Omega_2, \ldots$. We introduce the increments of $Y(t)$ between two consecutive inspections, i.e., $Z_n \coloneq Y(\bar\Omega_n) - Y(\bar\Omega_{n-1})$. Note that $(Z_n)_n$ is a sequence with i.i.d.\ random variables, say, distributed as the generic random variable $Z$. As a consequence, we can rewrite our target probability as 
\[p(u\,|\,\Omega)={\mathbb P}(T(u)<\infty) = {\mathbb P}(\tau(u)<\infty),\]
where
\begin{align*}
T(u) &\coloneqq  \inf \left\{n\in{\mathbb N}: \sum_{i = 1}^n Z_i > u\right\},\:\:\:\:
\tau(u) \coloneqq  \inf_{n\in{\mathbb N}:Y(\bar\Omega_n) > u } \bar\Omega_n,
\end{align*}
so that $Y(\tau(u)) = \sum_{i=1}^{T(u)}Z_i.$ In this section we assume that the {\it net profit condition} is in place, i.e., ${\mathbb E}\,Z<0$ or equivalently $\varphi'(0)>0$, so that bankruptcy becomes increasingly rare as $u\to\infty.$

To find a convenient expression for $p(u\,|\,\Omega)$, we work with an alternative measure $\QQ$. 
We assume that we are in the {\it light-tailed regime}, in that there is a $\theta^\star>0$ such that
${\mathbb E}\,e^{\theta^\star Z}=1$. Then the probability measure $\QQ$ is defined through the relation
\begin{equation}\label{QZ}\QQ(Z\in {\rm d}x) = \p(Z\in {\rm d}x) \,e^{\theta^\star x},\end{equation}
which integrates to 1 by the definition of $\theta^\star$. Notice that
\begin{align*} {\mathbb E}\,e^{\theta Z}&= \int_{0}^\infty \p(\Omega\in {\rm d}t) \,{\mathbb E}\,e^{\theta Y(t)}= \int_{0}^\infty \p(\Omega\in {\rm d}t) \,e^{\varphi(-\theta)\,t}= {\mathbb E}\,e^{\varphi(-\theta)\,\Omega}.
\end{align*}
As $\Omega\geqslant 0$, to obtain that ${\mathbb E}\,e^{\theta^\star Z}=1$ one should necessarily have that $\varphi(-\theta^\star) = 0$, implying that $\QQ$ is {\it not} affected by the distribution of the inter-inspection times $\Omega.$ Note that due to the convexity of $\varphi(\cdot)$ and $\varphi'(0)>0$, we have that $\varphi'(-\theta^\star)<0$; we throughout assume that $|\varphi'(-\theta^\star)|<\infty.$ Under the measure ${\mathbb Q}$, the Laplace exponent of the driving L\'evy process becomes $\varphi_{\mathbb Q}(\alpha)\coloneqq \varphi(\alpha-\theta^\star)$; in the sequel $\psi_{\mathbb Q}(\cdot)$ denotes the right-inverse of $\varphi_{\mathbb Q}(\cdot)$. 

Let $\ell(x)$ denote the likelihood ratio ${\rm d}\p/{\rm d}{\mathbb Q}$ evaluated in the argument $x$. Then we have the fundamental identity \cite[Section XIII.3]{ASM2}, describing the `translation' from sampling under the alternative measure $\QQ$ to sampling under the actual measure $\p$:
\[{\mathbb P}(T(u)<\infty) = {\mathbb E}\,1\{T(u)<\infty\} = {\mathbb E}_\QQ\big[1\{T(u)<\infty\}\,\ell(Z_1)\cdots\ell(Z_{T(u)})\big],\]
where we also use the independence of the random variables $Z_i$.

It is readily verified that
\[\frac{\rm d}{{\rm d}\theta}{\mathbb E}_\QQ e^{\theta Z} = -\varphi'(-\theta-\theta^\star)\,{\mathbb E}\,\big[\Omega \,e^{\varphi(-\theta-\theta^\star)\Omega}\big],\]
so that
\[{\mathbb E}_\QQ Z = \left.\frac{\rm d}{{\rm d}\theta}{\mathbb E}_\QQ \,e^{\theta Z}\right|_{\theta=0}=-\varphi'(-\theta^\star)\,{\mathbb E}\,\big[\Omega \,e^{\varphi(-\theta^\star)\Omega}\big] >0.\]
This means that $T(u)<\infty$ almost surely under ${\mathbb Q}$ (i.e., bankruptcy is not rare under ${\mathbb Q}$), entailing that 
\[{\mathbb P}(T(u)<\infty) =  {\mathbb E}_\QQ\big[\ell(Z_1)\cdots\ell(Z_{T(u)})\big].\]
Due to the fact that $\ell(Z_i) = e^{-\theta^\star Z_i}$, upon combining the above, we thus find that
\[p(u\,|\,\Omega) = {\mathbb E}_\QQ\big[\ell(Z_1)\cdots\ell(Z_{T(u)})\big] = {\mathbb E}_\QQ\exp\left(-\theta^\star \sum_{i=1}^{T(u)} Z_i\right) = {\mathbb E}_\QQ\,e^{-\theta^\star Y(\tau(u))}. \]
We have arrived at the following result.

\begin{theorem}\label{CoM}
For any $u\geqslant 0$,
\[p(u\,|\,\Omega) = {\mathbb E}_{\mathbb Q}\, e^{-\theta^\star Y(\tau(u))}.\]
\end{theorem}

A consequence of this result is the following uniform bound on $p(u\,|\,\Omega)$, which can be seen as a Lundberg-type inequality for our model with i.i.d.\ inter-inspection times. It follows from the obvious property $Y(\tau(u))> u$. Remarkably, this upper bound does not depend on the inter-inspection time distribution (and will therefore be a loose bound in case of infrequent inspections).

\begin{corollary} For any $u\geqslant 0$, irrespective of the inter-inspection time $\Omega$,
\[p(u\,|\,\Omega)\leqslant e^{-\theta^\star  u}.\]
\end{corollary}

\subsection{Exact asymptotics}
In the remainder of this section, our objective is to use Theorem \ref{CoM} to identify the exact asymptotics of $p(u\,|\,\Omega)$. More concretely, we aim to identify the constant $\gamma[\Omega]$ featuring in
\begin{equation}
    \label{exg} 
\lim_{u\to\infty}p(u\,|\,\Omega)\, e^{\theta^\star  u} \eqqcolon \gamma[\Omega],\end{equation}
for the inter-inspection times $\Omega$ studied in Sections \ref{Exp_HypExp} and \ref{Erl_HypErl}. In view of Theorem \ref{CoM}, this means that, with $R(u)$ denoting the {\it overshoot} $Y(\tau(u))-u$, we wish to evaluate the object
\[\gamma[\Omega] = \lim_{u\to\infty} {\mathbb E}_\QQ e^{-\theta^\star R(u)}.\]
Following \cite{KMD2}, the idea is to evaluate this quantity using that, with $U_\alpha$ denoting an exponentially distributed random variable with parameter $\alpha$ (i.e., with mean $\alpha^{-1}$, having probability density function $\alpha e^{-\alpha u}$), 
\begin{equation}\label{gom}\gamma[\Omega] = \lim_{\alpha\downarrow 0} {\mathbb E}_\QQ e^{-\theta^\star R(U_\alpha)}=
 \lim_{\alpha\downarrow 0}\alpha\,\pi^\circ(\alpha\,|\,\Omega),\end{equation}
 where the transform
 \[\pi^\circ(\alpha\,|\,\Omega) \coloneqq  \int_0^\infty e^{-\alpha u} \,{\mathbb E}_\QQ e^{-\theta^\star R(u)}\,{\rm d}u\]
 will be evaluated essentially relying on the techniques developed in Sections \ref{Exp_HypExp} and \ref{Erl_HypErl}.
 In this derivation we heavily use the `master equation', with $Z$ distributed as $Y(T_\omega)$,
 \begin{align}
 \pi^\circ(\alpha\,|\,\Omega) =\:& \int_0^\infty e^{-\alpha u} \int_u^\infty {\mathbb Q}(Z\in{\rm d}v)\,e^{-\theta^\star(v-u)}\,{\rm d}u\:+\notag\\
 &\int_0^\infty e^{-\alpha u} \int_{-\infty}^u {\mathbb Q}(Z\in{\rm d}v)\,{\mathbb E}_\QQ e^{-\theta^\star R(u-v)}\,{\rm d}u,\label{master}
 \end{align}
 where the first term corresponds to the scenario of the first $Z$ directly exceeding level $u$, and the second term to the scenario that the first $Z$ is still below $u$; cf.\ Equation \eqref{eq1}.
 \subsubsection*{-- Exponentially distributed inter-inspection times}
 Our first goal is to evaluate $\pi^\circ(\alpha\,|\,\Omega)$ for exponentially distributed inter-inspection times with parameter $\omega$, essentially following the approach that we developed in Section~\ref{Exp inter-inspection}. The underlying idea is to separate $Z$ into the independent parts $Z^+\coloneqq \bar Y(T_\omega)$ and $Z^-\coloneqq -\underline Y(T_\omega)$, the latter random variable being exponentially distributed with parameter $\theta\coloneqq \psi_{\mathbb Q}(\omega)$. As a consequence, the first term in the right hand side of \eqref{master} can be rewritten as
 \begin{align*}
\int_0^\infty e^{-\alpha u} &\int_u^\infty {\mathbb Q}(Z\in{\rm d}v)\,e^{-\theta^\star(v-u)}\,{\rm d}u\\
&=\int_{u=0}^\infty e^{-\alpha u} \int_{v=u}^\infty \int_{w=0}^\infty \theta e^{-\theta w}\,{\mathbb Q}(Z^+-w\in{\rm d}v)\,e^{-\theta^*(v-u)}\,{\rm d}w\,{\rm d}u\\
&=\int_{u=0}^\infty e^{-\alpha u} \int_{v=u}^\infty \int_{x=v}^\infty \theta e^{-\theta (x-v)}\,{\mathbb Q}(Z^+\in{\rm d}x)\,e^{-\theta^*(v-u)}\,{\rm d}v\,{\rm d}u,
 \end{align*}
 where in the last equality the change of variable $x\coloneqq v+w$ has been performed. Swapping the inner two integrals and computing the integral over $v$ yields that the expression in the previous display reads
 \[\int_{u=0}^\infty e^{-\alpha u} \int_{x=u}^\infty \frac{\theta}{\theta^\star-\theta}\big(e^{-\theta(x-u)}-e^{-\theta^\star(x-u)}\big) \,{\mathbb Q}(Z^+\in{\rm d}x)\diff u.\]
 After interchanging the two remaining integrals and integrating over $u$, we readily obtain
 \[\frac{\theta}{\theta^\star-\theta}\int_0^\infty\left(\frac{e^{-\theta x}-e^{-\alpha x}}{\alpha-\theta}-\frac{e^{-\theta^\star x}-e^{-\alpha x}}{\alpha-\theta^\star}\right)\,{\mathbb Q}(Z^+\in{\rm d}x).\]
We thus end up with the first integral in the right hand side of \eqref{master} equalling
\[\int_0^\infty e^{-\alpha u} \int_u^\infty {\mathbb Q}(Z\in{\rm d}v)\,e^{-\theta^\star(v-u)}\,{\rm d}u=F^\circ(\alpha,\omega),\]
where, recalling that $\theta=\psi_{\mathbb Q}(\omega)$,
\begin{equation}\label{Fcirc}F^\circ(\alpha,\omega)\coloneqq \frac{\theta}{\theta^\star-\theta}\left(\frac{\xi_{\mathbb Q}(\theta,\omega)-\xi_{\mathbb Q}(\alpha,\omega)}{\alpha-\theta}-
\frac{\xi_{\mathbb Q}(\theta^\star,\omega)-\xi_{\mathbb Q}(\alpha,\omega)}{\alpha-\theta^\star}\right)\end{equation}
and, using that $\psi_{\mathbb Q}(\beta)=\psi(\beta)+\theta^\star$ for any $\beta\geqslant 0$,
\begin{align*}\xi_{\mathbb Q}(\alpha,\beta)&\coloneqq {\mathbb E}_{\mathbb Q}e^{-\alpha \bar Y(T_\beta)}= {\mathbb E}_{\mathbb Q}\,e^{-\alpha Z^+}=\frac{\alpha -\psi_{\mathbb Q}(\beta)}{\varphi_{\mathbb Q}(\alpha)-\beta}\frac{\beta}{\psi_{\mathbb Q}(\beta)}=\frac{\psi(\beta)+\theta^\star-\alpha}{\psi(\beta)+\theta^\star}\frac{\beta}{\beta-\varphi(\alpha-\theta^\star)};\end{align*}
cf.\ \eqref{xidef}. Regarding the second integral in the right hand side of \eqref{master}, the procedure step-by-step mimics the one developed in Section \ref{Exp inter-inspection}, thus leading to
    \[\frac{\theta}{\theta-\alpha}\ \big({\pi}^\circ(\alpha\,|\,T_\omega)\,\xi_{\mathbb Q}(\alpha,\omega)
    - {\pi}^\circ(\theta\,|\,T_\omega)\,\xi_{\mathbb Q}(\theta,\omega)\big).\]

Combining the above, we obtain the equation
\[{\pi}^\circ(\alpha\,|\,T_\omega)\left(1-
\frac{\theta}{\theta-\alpha}\,\xi_{\mathbb Q}(\alpha,\omega)
\right)=F^\circ(\alpha,\omega)-G^\circ(\alpha,\omega),\:\:\:\:
G^\circ(\alpha,\omega)
\coloneqq
{\theta}\,\xi_{\mathbb Q}(\theta,\omega)\frac
{{\pi}^\circ(\theta\,|\,T_\omega)}{\theta-\alpha}.\]
We are left with determining the unknown constant ${\pi}^\circ(\theta\,|\,T_\omega)$. To this end, it is first observed that  $\alpha = \theta^\star$ is a root of the left hand side, and hence it is a root of the right hand side as well. This leads to
\[{\pi}^\circ(\theta\,|\,T_\omega) =  \frac{\theta-\theta^\star}{\theta}\frac{F^\circ(\theta^\star,\omega)}{\xi_{\mathbb Q}(\theta,\omega)}.\]
We conclude that
\begin{equation}
    \label{picirc}
{\pi}^\circ(\alpha\,|\,T_\omega) = \frac{(\theta-\alpha)F^\circ(\alpha,\omega)-\theta\,{\pi}^\circ(\theta\,|\,T_\omega)\,\xi_{\mathbb Q}(\theta,\omega) }{\theta-\alpha-\theta\,\xi_{\mathbb Q}(\alpha,\omega)}.\end{equation}
Now that we have found the transform ${\pi}^\circ(\alpha\,|\,T_\omega)$, we are in the position to identify the constant $\gamma[{T_\omega}]$. Indeed, appealing to relation \eqref{gom}, we find as a direct application of L'H\^opital's rule, that the expression for ${\pi}^\circ(\alpha\,|\,T_\omega)$ displayed in \eqref{picirc} leads to
\begin{align*}\gamma[{T_\omega}] &= \lim_{\alpha\downarrow 0} \alpha \,{\pi}^\circ(\alpha\,|\,T_\omega)= \lim_{\alpha\downarrow 0} \frac{\alpha\,\kappa}{{\theta-\alpha-\theta\,\xi_{\mathbb Q}(\alpha,\omega)}}=\frac{\kappa}{-\big(1+\theta f(0,\omega)\big)},
\end{align*}
where
\begin{align*}f(0,\omega)&\coloneqq  \left.\frac{\rm d}{{\rm d}\alpha} \xi_{\mathbb Q}(\alpha,\omega)\right|_{\alpha= 0}=\frac{1}{\theta}\frac{\theta\,\varphi'_{\mathbb Q}(0)-\omega}{\omega},\\
\kappa&\coloneqq  \theta\big( F^\circ(0,\omega) -\pi^\circ(\theta\,|\,T_\omega)\,\xi_{\mathbb Q}(\theta,\omega)\big)= \theta\,F^\circ(0,\omega)+(\theta^\star-\theta)\,F^\circ(\theta^\star,\omega).\end{align*}
Observe that, using that $\theta=\psi_{\mathbb Q}(\omega)=\psi(\omega)+\theta^\star$,
\[1+\theta f(0,\omega) =\frac{\theta\,\varphi_{\mathbb Q}'(0)}{\omega}. \]
Our next task is to write the constant $\kappa$ in a more explicit form. By computing $F^\circ(0,\omega)$ and $F^\circ(\theta^\star,\omega)$, it is directly seen that
\begin{align*}\kappa =\frac{\theta}{\psi(\omega)}\Big(&\xi_{\mathbb Q}(\theta,\omega)-1-\frac{\theta}{\theta^\star}\left(\xi_{\mathbb Q}(\theta^\star,\omega)-1\right)-\xi_{\mathbb Q}(\theta,\omega)+\xi_{\mathbb Q}(\theta^\star,\omega)+\psi(\omega)\,f(\theta^\star,\omega)\Big),\end{align*}
where
\[\xi_{\mathbb Q}(\theta^\star,\omega)=\frac{\psi(\omega)}{\theta},\:\:\:f(\theta^\star,\omega)\coloneqq  \left.\frac{\rm d}{{\rm d}\alpha} \xi_{\mathbb Q}(\alpha,\omega)\right|_{\alpha =\theta^\star}=\frac{1}{\theta}\frac{\psi(\omega)\,\varphi'(0)-\omega}{\omega}.\]
Putting the above results together, standard algebra yields that $\kappa=-\psi(\omega)\,\varphi'(0)/\omega$. This leads us to the following conclusion, in line with \cite[Proposition 4.1]{BM2}. 
\begin{theorem} \label{thEXP} As $u\to\infty$, we have that $p(u\,|\,T_\omega)\, e^{\theta^\star  u}\to \gamma[{T_\omega}]$, where
\[\gamma[{T_\omega}] =-\frac{\varphi'(0)}{\varphi'_{\mathbb Q}(0)}\frac{\psi(\omega)}{\theta}=  -\frac{\varphi'(0)}{\varphi'_{\mathbb Q}(0)}\frac{\psi(\omega)}{\psi(\omega)+\theta^\star}.\]
\end{theorem}
As the inspection rate $\omega$ grows large, we recover the asymptotics under permanent inspection. Indeed, in agreement with the findings of \cite{DON}, revealing the exact asymptotics for spectrally-positive L\'evy processes, we obtain that $\gamma[{T_\omega}]\to -{\varphi'(0)}/{\varphi'_{\mathbb Q}(0)}$ as $\omega\to\infty$.

 \subsubsection*{-- Hyperexponentially distributed inter-inspection times}
Denoting as before $\theta_i\coloneqq \psi_{\mathbb Q}(\omega_i)$, following the same approach as in the case of exponentially distributed inter-inspection times, we obtain the relation
\begin{align*}{\pi}^\circ(\alpha\,|\,T_{{\bs \omega},{\bs p}})&\left(1-\sum_{i=1}^d p_i
\frac{\theta_i}{\theta_i-\alpha}\,\xi_{{\mathbb Q}}(\alpha,\omega_i)
\right)=\sum_{i=1}^d p_i\left(F^\circ(\alpha,\omega_i)-
G^\circ(\alpha,\omega_i)
\right),\end{align*}
with $F^\circ(\cdot,\cdot)$ as defined in \eqref{Fcirc} with $\theta$ replaced by $\theta_i$,
and
\[G^\circ(\alpha,\omega_i)\coloneqq{\theta_i}\,\xi_{{\mathbb Q}}(\theta_i,\omega_i)
\frac{{\pi}^\circ(\theta_i\,|\,T_{{\bs \omega},{\bs p}})}{\theta_i-\alpha}.\]
Inserting the expression for $\xi_{{\mathbb Q}}(\alpha,\omega_i)$, we thus obtain 
\begin{equation}\label{hypexppicirc}
    {\pi}^\circ(\alpha\,|\,T_{{\bs \omega},{\bs p}})=\frac{\displaystyle \sum_{i=1}^d p_i\left(F^\circ(\alpha,\omega_i)-
\frac{\theta_i}{\theta_i-\alpha}
\,\xi_{{\mathbb Q}}(\theta_i,\omega_i)\,z_i\right)}{\displaystyle 1-\sum_{i=1}^d p_i
\frac{\omega_i}{\omega_i-\varphi(\alpha-\theta^\star)}}.
\end{equation}
where the $d$ unknowns $z_i\coloneqq {\pi}^\circ(\theta_i\,|\,T_{{\bs \omega},{\bs p}})$ are yet to be determined. The denominator in the right-hand side of the previous display (and hence the numerator as well) has precisely $d$ roots in the right half of the complex plane, say $\alpha_1,\ldots,\alpha_d$. This follows by the same line of reasoning as in the proof of Lemma \ref{rootshypexp}. We consider the same MAP, but now we have $\varphi_1(\alpha)=\dots=\varphi_d(\alpha)=\varphi_\mathbb{Q}(\alpha)$ and $\boldsymbol{\eta}=0$ (i.e., no killing). From \cite[Theorem 2]{IBM} we know that, since we have positive drift under $\mathbb{Q}$, the equation has exactly $d$ roots in the right half of the complex plane. Applying the same reasoning as in Remark \ref{R1}, the roots $\alpha_1,\ldots,\alpha_d$ are real-valued, and hence the solution ${\bs z}$ of the linear system, with $j=1,\ldots,d$,
\[\sum_{i=1}^d p_iF^\circ(\alpha_j,\omega_i)=\sum_{i=1}^d p_i
\frac{\theta_i}{\theta_i-\alpha_j}
\,\xi_{{\mathbb Q}}(\theta_i,\omega_i)\,z_i\]
is real-valued as well. Defining 
\[\Phi(\alpha)\coloneqq \prod_{i=1}^d (\omega_i -\varphi(\alpha-\theta^\star)),\:\:\:\:
\Phi_i(\alpha)\coloneqq \frac{\Phi(\alpha)}{\omega_i-\varphi(\alpha-\theta^\star)},\]
and applying the same reasoning as in the exponential case,
we have that the expression for $\pi^{\circ}(\alpha\,|\,T_{\omega,\bs{p}})$ displayed in \eqref{hypexppicirc} leads to
\[
\gamma[T_{\omega,\bs{p}}]= \lim_{\alpha\downarrow 0} \alpha \,{\pi}^\circ(\alpha\,|\,T_{\omega,\bs{p}})= \lim_{\alpha\downarrow 0}\frac{\alpha K(\alpha)\,\Phi(\alpha)}{\Phi(\alpha)-\sum_{i=1}^d p_i \omega_i \Phi_i(\alpha)},
\]
where 
\[
K(\alpha)\coloneqq \sum_{i=1}^d p_i\big(F^\circ(\alpha,\omega_i) -\frac{\theta_i}{\theta_i-\alpha}\xi_{\mathbb Q}(\theta_i ,\omega_i)\,z_i\big).
\]
We obtain the following result.
\begin{theorem} \label{asshe} As $u\to\infty$, we have that $p(u\,|\,T_{{\bs \omega},{\bs p}})\, e^{\theta^\star  u}\to \gamma[{T_{{\bs \omega},{\bs p}}}]$, where 
\[\gamma[{T_{{\bs \omega},{\bs p}}}] =\frac{\displaystyle \Phi(0)\sum_{i=1}^d p_i\big(F^\circ(0,\omega_i) -\xi_{\mathbb Q}(\theta_i ,\omega_i)\,z_i\big)}{\displaystyle \Phi'(0)-\sum_{i=1}^dp_i\Phi_i'(0)\,\omega_i}. \]
\end{theorem}

 \subsubsection*{-- Erlang distributed inter-inspection times}
 In the Erlang case the computations become somewhat more involved. It is readily verified that the first term in the right hand side of \eqref{master} can be rewritten as
 \begin{align}
\int_{u=0}^\infty e^{-\alpha u} \int_{x=u}^\infty\int_{v=u}^x   e^{-\theta (x-v)}\,\frac{\theta^k{(x-v)^{k-1}}}{(k-1)!}\,{\mathbb Q}(Z_k^+\in{\rm d}x)\,e^{-\theta^*(v-u)}\,{\rm d}v\,{\rm d}u.\label{erlpicirc}
 \end{align}
Then we use that
 \[\int_u^x e^{-(\theta^\star-\theta)v}\frac{\theta^k (x-v)^{k-1}}{(k-1)!}{\rm d}v = \left(\frac{\theta}{\theta-\theta^\star}\right)^k\left(e^{-(\theta^\star-\theta)x}-\sum_{n=0}^{k-1}e^{-(\theta^\star-\theta)u}\frac{((\theta^\star-\theta)(x-u))^n}{n!}\right),\]
so that \eqref{erlpicirc} reads, after changing the order of integration,
\begin{align*}
\int_{x=0}^\infty \int_{u=0}^x &\left(e^{-\theta^\star x}e^{-(\alpha-\theta^\star)u}{\rm d}u-\sum_{n=0}^{k-1}\int_{u=0}^xe^{-\alpha x}e^{-(\theta-\alpha)u}\frac{((\theta^\star-\theta)u)^n}{n!}{\rm d}u\right)\left(\frac{\theta}{\theta-\theta^\star}\right)^k {\mathbb Q}(Z_k^+\in{\rm d}x).
\end{align*}
After straightforward calculus we conclude that the first term in the right-hand side of \eqref{master} equals, with $\xi_{\mathbb Q,k}(\alpha,\beta)$ defined in the evident manner,
\begin{align*}
I^\circ(\alpha,\omega)&\coloneqq\left(\frac{\theta}{\theta-\theta^\star}\right)^k\frac{\xi_{{\mathbb Q},k}(\theta^\star,\omega)-\xi_{{\mathbb Q},k}(\alpha,\omega) }{\alpha-\theta^\star}\:-\\
&\:\:\:\:\left(\frac{\theta}{\theta-\theta^\star}\right)^k\sum_{n=0}^{k-1}\frac{(\theta-\theta^\star)^n}{(\theta-\alpha)^{n+1}}\left(\xi_{\mathbb Q,k}(\alpha,\omega)-\sum_{m=0}^n\xi_{\mathbb Q,k}^{(m)}(\theta,\omega)\frac{(\alpha-\theta)^m}{m!}\right).
\end{align*}
The second term in the right hand side of \eqref{master}
can be found analogously to Lemma \ref{LL3}. Defining, recalling the definition of $\bar J_i(\alpha,\beta,\omega)$ from \eqref{defjbar}, 
\begin{align*}\bar J_i^\circ(\alpha,\omega)&\coloneqq \sum_{n=i}^{k-1}\delta^{\mathbb Q}_{k-1-n,k}\frac{(\alpha-\theta)^i}{i!}\left(\frac{\theta}{\theta-\alpha}\right)^{n+1},\:\:\:
\delta_{n,k}^{\mathbb Q}\coloneqq \frac{(-\theta)^n}{n!}\,\xi_{\mathbb Q,k}^{(n)}(\theta,\omega),\end{align*}
we end up with, writing $J^{\circ}(\alpha,\omega)\coloneqq\sum_{i=0}^{k-1}\bar{J}_i^{\circ}(\alpha,\omega)\,z_i$ with $z_i\coloneqq (\pi^\circ)^{(i)}(\theta\,|\,E_{k,\omega})$,
\[\pi^\circ(\alpha\,|\,E_{k,\omega}) = I^\circ(\alpha,\omega) - J^{\circ}(\alpha,\omega)+\left(\frac{\omega}{\omega-\varphi(\alpha-\theta^\star)}\right)^k\pi^\circ(\alpha\,|\,E_{k,\omega}),\]
or equivalently 
\[\pi^\circ(\alpha\,|\,E_{k,\omega})=\frac{\displaystyle I^\circ(\alpha,\omega) -  J^\circ(\alpha,\omega)}{\displaystyle 1- \left(\frac{\omega}{\omega-\varphi(\alpha-\theta^\star)}\right)^k},\]
where we still need to identify, for $i=0,\ldots,k-1$, the constants
$z_i$.

As before, an application of \cite[Theorem 2]{IBM} reveals that the denominator has precisely $k$ roots in the right half of the complex plane, say $\alpha_1,\ldots,\alpha_k,$ and hence the corresponding numerator as well. It implies that ${\bs z}$ can be determined by solving, for $j=1,\ldots,k$,
\[I^\circ(\alpha_j,\omega) = \sum_{i=0}^{k-1}\bar J_i^\circ(\alpha_j,\omega)\,z_i.\]
We obtain the following result, by following the same procedure as before; observe that $\varphi'(-\theta^\star)<0.$

\begin{theorem} \label{asserl} As $u\to\infty$, we have that $p(u\,|\,E_{k,\omega})\, e^{\theta^\star  u}\to \gamma[{E_{k,\omega}}]$, where
\[\gamma[E_{k,\omega}] =\frac{\displaystyle I^\circ(0,\omega)-  J^\circ(0,\omega)}{-k\,\varphi'(-\theta^\star)/\omega}. \]
\end{theorem}

  \subsubsection*{-- Hyper-Erlang distributed inter-inspection times}
Following the same approach as in the case of Erlang distributed inter-inspection times, we obtain the relation
\begin{align*}
    \pi^\circ(\alpha\,|\,E_{\bs{k,p},\omega}) = \frac{\displaystyle \sum_{i=1}^d p_i\big( I_i^\circ(\alpha,\omega) -  J_{i}^\circ(\alpha,\omega)\big)}{\displaystyle 1-\sum_{i=1}^d p_i\left(\frac{\omega}{\omega-\varphi(\alpha-\theta^\star)}\right)^{k_{i}}}
    ,
\end{align*}
where
\begin{align*}
    I_i^\circ(\alpha,\omega)&\coloneqq \left(\frac{\theta}{\theta-\theta^\star}\right)^{k_i}\frac{\xi_{{\mathbb Q},k_i}(\theta^\star,\omega)-\xi_{{\mathbb Q},k_i}(\alpha,\omega) }{\alpha-\theta^\star}\:-\\
&\:\:\:\:\left(\frac{\theta}{\theta-\theta^\star}\right)^{k_i}\sum_{n=0}^{k_i-1}\frac{(\theta-\theta^\star)^n}{(\theta-\alpha)^{n+1}}\left(\xi_{\mathbb Q,k_i}(\alpha,\omega)-\sum_{m=0}^n\xi_{\mathbb Q,k_i}^{(m)}(\theta,\omega)\frac{(\alpha-\theta)^m}{m!}\right),
\end{align*}
and, with $z_j\coloneqq (\pi^\circ)^{(j)}(\theta\,|\,E_{\bs{k,p},\omega})$,
\begin{align*}
    \bar{J}^\circ_{i,j}(\alpha,\omega)\coloneqq \sum_{n=j}^{k_i-1}\delta^{\mathbb{Q}}_{k_i-1-n,k_i}\frac{(\alpha-\theta)^j}{j!}\left(\frac{\theta}{\theta-\alpha}\right)^{n+1},\:\:\:
    J_i^\circ(\alpha,\omega)\coloneqq \sum_{j=0}^{k_i-1}\bar{J}^\circ_{i,j}(\alpha, \omega)\,z_j.
\end{align*}
As before assuming (without loss of generality) that $k_1<\cdots<k_d$, we are left with the task to identify the constants $z_j$, for $j=0,\dots,k_d-1$. Again, application of \cite[Theorem 2]{IBM} reveals that the denominator has exactly $k_d$ roots in the right half of the complex plane, say $\alpha_1,\dots,\alpha_{k_d}$, which are then roots of the numerator as well. Hence, the vector  ${\bs z}$ can be determined by solving the linear system
\begin{align*}
    \sum_{i=1}^d p_iI_i^\circ (\alpha_\ell, \omega)= \sum_{i=1}^d p_i\sum_{j=0}^{k_i-1}\bar{J}^\circ_{i,j}(\alpha_\ell, \omega)\,z_j,
\end{align*}
for  $\ell=1,\dots,k_d$.
We obtain the following result.
\begin{theorem}\label{asshyperl}
    As $u\to\infty$, we have that $p(u\,|\,E_{\bs{k,p},\omega})\, e^{\theta^\star  u}\to \gamma[{E_{\bs{k,p},\omega}}]$, where
\[
\gamma[{E_{\bs{k,p},\omega}}]=\frac{\displaystyle \sum_{i=1}^d p_i\big( I_i^\circ(0,\omega)-J_i^\circ(0,\omega)\big)}{
\displaystyle -\sum_{i=1}^d p_ik_i \,\varphi'(-\theta^\star)/\omega}.
\]
\end{theorem}

\section{Numerical evaluation, rare-event simulation}\label{RES}
The objective of this section is to demonstrate how to numerically evaluate $p(u\,|\,\Omega)$, the all-time ruin probability at inspection epochs. 
Two techniques are highlighted: a highly efficient importance sampling technique, and the numerical evaluation of the Cram\'er-Lundberg type asymptotics.
We start by explaining our importance sampling procedure, and prove that it is logarithmically efficient. Then we introduce the class of instances considered in our experiments, and point out how to numerically evaluate the constants $\theta^\star$ and $\gamma[\Omega]$. We conclude this section by presenting the output of our experiments.

\subsection{Importance sampling}
As we observed before, the process $(Y(\bar\Omega_n))_{n\in{\mathbb N}}$ constitutes a random walk with i.i.d.\ increments distributed as the random variable $Z$. The results of \cite{KOR} reveal that in the light-tailed regime introduced in Section \ref{CLA}, we have that \eqref{CLA} holds, suggesting the approximation 
\begin{equation}
    \label{APPR}
p(u\,|\,\Omega)\approx \gamma[\Omega]\, e^{-\theta^\star u}.\end{equation}
This approach, however, has two intrinsic difficulties. In the first place, the exact asymptotics relate to the regime that $u$ grows large, so that the approximation is expected to be accurate for large $u$. Hence, in the absence of error bounds, it is not a priori clear whether for the value of $u$ under consideration,  \eqref{APPR} provides an accurate approximation. In the second place, while we have succeeded to compute $\gamma[\Omega]$ for various types of inter-inspection times $\Omega$, no expression is available for general $\Omega$. As detailed in \cite[Theorem 3]{KOR}, 
with
\[T\coloneqq \inf\left\{n\in{\mathbb N}:\sum_{i=1}^n Z_i>0\right\},\:\:\:\:\:a\coloneqq {\mathbb E}\big[Z\,e^{\theta^\star Z}1\{T<\infty\}\big],\]
we do have that
\[\gamma[\Omega] = \frac{1-p(0\,|\,\Omega)}{\theta^\star\,a},\]
but it is not clear how to compute, for an arbitrarily distributed inter-inspection time $\Omega$, the quantities $p(0\,|\,\Omega)$ and $a$.

The above considerations show that there is a need to develop computational tools to efficiently and accurately evaluate $p(u\,|\,\Omega)$ for any non-negative random variable $\Omega$ and $u>0$. The numerical procedure we discuss in this subsection is based on simulation: we devise an importance sampling algorithm that applies the change of measure that we worked with in Section \ref{CLA}.
It succeeds in efficiently generating precise estimates for $p(u\,|\,\Omega)$, even in the  regime that $u$ grows large (and ruin is rare). 

The starting point of the algorithm is Theorem \ref{CoM}. In each (independent) run, the model is simulated under ${\mathbb Q}$ until the hitting time $\tau(u)$, where we recall that $u$ is reached almost surely. Let $W_i$ denote the value of $Y(\tau(u))$ in the $i$-th run. Our Siegmund-type estimator \cite[Section VI.2a]{AG}, based on $n$ runs, is
\begin{equation}\label{esti}p_n(u) \coloneqq \frac{1}{n} \sum_{i=1}^n e^{-\theta^\star W_i},\end{equation}
where the sample variance of this estimator is
\[v_n(u) \coloneqq \frac{1}{n-1}\left(\sum_{i=1}^n (e^{-\theta^\star W_i}-p_n(u))^2\right).\]
Noting that $W_i\geqslant u$, from
\[{\mathbb V}{\rm ar}_{\mathbb Q}\, e^{-\theta^\star W_i} = {\mathbb E}_{\mathbb Q} e^{-2\theta^\star W_i}- ({\mathbb E}_{\mathbb Q} e^{-\theta^\star W_i})^2\leqslant {\mathbb E}_{\mathbb Q} e^{-2\theta^\star W_i} \leqslant e^{-2\theta^\star u} \]
the following result follows directly; see \cite[Section VI.1]{AG} for the definition of logarithmic efficiency. 

\begin{theorem}
    The estimator \eqref{esti} is logarithmically efficient as $u\to\infty.$
\end{theorem}

\subsection{Change of measure}
We proceed by discussing how to sample the process $Y(\cdot)$ according to the alternative measure ${\mathbb Q}.$ As was pointed out in Section \ref{CLA}, under ${\mathbb Q}$ the distribution of $Z$ should follow \eqref{QZ}. We have that
\begin{align*}
    \e[e^{- \alpha Z}] &= \int_0^{\infty}  e^{\varphi(\alpha)\,t}\,{\mathbb P}(\Omega\in{\rm d}{t}) = {\mathbb E}\,e^{\varphi(\alpha)\,\Omega},
\end{align*}
so that
\[\e_{\mathbb Q}[e^{- \alpha Z}]  = {\mathbb E}\,e^{\varphi(\alpha-\theta^\star)\,\Omega}.\]
We conclude that under this new measure the net cumulative claim  process $Y(\cdot)$ has Laplace exponent $\varphi(\alpha - \theta^\star )$ instead of $\varphi(\alpha)$, while the distribution of $\Omega$ remains unaltered. This leaves us with the question of how to sample a L\'evy process with Laplace exponent $\varphi(\alpha - \theta^\star )$. We consider the practically relevant case of $Y(\cdot)$ corresponding to the sum of two independent processes: a compound Poisson process $C(\cdot)$ and a Brownian motion $B(\cdot)$; it is noted that  any spectrally positive L\'evy process can be approximated arbitrarily closely by such a sum \cite{AR}. 

Let the Brownian motion have variance coefficient $\sigma^2$, and let the compound Poisson process have drift $r$, arrival rate $\lambda$, with positive jump sizes characterized by the Laplace-Stieltjes transform $b(\cdot).$ This means that
\[\varphi(\alpha) = \log {\mathbb E}\exp(-\alpha B(1)-\alpha C(1)) = \frac{\alpha^2\sigma^2}{2}+r\alpha-\lambda(1-b(\alpha)), \]
where we assume that the net profit condition is in place (i.e., ${\mathbb E}\,Y(1)<0$, or equivalently $-\lambda \,b'(0)<r$). To identify $Y(\cdot)$ under ${\mathbb Q}$, observe that
\begin{align*}\varphi(\alpha - \theta^\star )=\varphi_{\mathbb Q}(\alpha)
&= \frac{(\alpha- \theta^\star)^2\sigma^2}{2}+r(\alpha- \theta^\star)-\lambda(1-b(\alpha- \theta^\star)) + \varphi(-\theta^\star)\\
&= \frac{\alpha^2\sigma^2}{2}+(r-\theta^\star\sigma^2)\alpha-\lambda(b(-\theta^\star)-b(\alpha-\theta^\star))\\
&= \frac{\alpha^2\sigma^2}{2}+(r-\theta^\star\sigma^2)\alpha-\lambda b(-\theta^\star)\left(1-\frac{b(\alpha-\theta^\star)}{b(-\theta^\star)}\right),
\end{align*}
recalling that $\varphi(-\theta^\star)=0.$ This means that under ${\mathbb Q}$ the Brownian part remains unchanged, but that the compound Poisson process now has drift $r-\theta^\star\sigma^2$, the arrival rate becomes $\lambda b(-\theta^\star)$, and the claim sizes are now characterized by the Laplace-Stieltjes transform $b(\cdot-\theta^\star)/b(-\theta^\star).$

In the experiments we let the claim sizes be exponentially distribution with parameter $\mu$. It is readily verified that (for $\alpha>-\mu$)
\[\varphi(\alpha) = \alpha\,\frac{\sigma^2\alpha^2 +(\mu\sigma^2 +2r)\alpha+2(r\mu-\lambda)}{2(\mu+\alpha)},\]
so that
\[\theta^\star = \frac{\mu\sigma^2 +2r}{2\sigma^2} - \frac{\sqrt{(\mu\sigma^2+2r)^2-8\sigma^2(r\mu-\lambda)}}{2\sigma^2}.\]
Under ${\mathbb Q}$ the claim sizes are exponentially distributed with parameter $\mu-\theta^\star.$ 

\subsection{Numerical experiments: two-moment approximation for inter-inspection times}
In this series of numerical experiments we assess $p(u\,|\,\Omega)$ for an instance in which the inter-inspection time $\Omega$ does {\it not} have one of the phase-type distributions covered in this paper, and compare it to its counterpart in which $\Omega$ is approximated by $\Omega_{\rm app}$ determined by the two-moment fit that was discussed in Section \ref{TMF}. 
For ease we use (without losing any generality) the normalization ${\mathbb E}\,\Omega = 1$, and vary ${\mathbb V}{\rm ar}\,\Omega$ (where we cover values both below and above 1). 
The inter-inspection time $\Omega$ we consider has a lognormal distribution,  with its two parameters chosen such that the mean is 1 and the variance has the target value.

We proceed by pointing out how to numerically evaluate the constant $\gamma[\Omega_{\rm app}]$ for these~$\Omega_{\rm app}.$
\begin{itemize}
\item[$\circ$]
If the coefficient of variation of the inter-inspection times is larger than 1, then the two-moment fit prescribes that we should use a hyperexponential random variable $T_{{\bs p},{\bs \omega}}$ with $d=2$, where the parameter vectors ${\bs p}$ and ${\bs \omega}$ are chosen such that the first two moments have the desired values (and in addition the so-called balanced means condition is fulfilled; see e.g.\ \cite[Section 3.1]{KUIP}). It requires some elementary calculus to verify that $\alpha_1$ and $\alpha_2$ are, respectively, $\theta^\star$ and the positive root of (with $p:=p_1$)
\[\varphi_{\mathbb Q}(\alpha) = (1-p)\omega_1+p\omega_2.\]
As is readily checked from the expression for $\varphi_{\mathbb Q}(\alpha)$, finding $\alpha_2$ requires solving an equation involving a third degree polynomial.
\item[$\circ$] If the coefficient of variation is larger than 1, then the two-moment fit entails that we use an $E_{{\bs k},{\bs p},\omega}$ distribution with $d=2$, $k_1=k$, $k_2=k+1$ and $p_1=p$, where $p$ and $k$ are such that the first two moments have the desired value. In this case $\alpha_1,\ldots,\alpha_{k+1}$ are the roots in the right half of the complex plane of the polynomial equation
\begin{equation}\label{rooteq}1-\left(\frac{\omega-\varphi_{\mathbb Q}(\alpha)}{\omega}\right)^{k+1}=\frac{p\varphi_{\mathbb Q}(\alpha)}{\omega},\end{equation}
where, evidently, $\theta^\star$ is one such root. There are two ways of finding the roots $\alpha_1,\ldots,\alpha_{k+1}$, both relying on the availability of a procedure to find the, potentially complex, roots of a polynomial equation, which is available in standard numerical packages. An approach to do so would be to insert the expression for $\varphi_{\mathbb Q}(\alpha)$ into \eqref{rooteq}, and to directly solve the resulting polynomial of order $3k+3$, after which one selects the roots with a positive real part. There is an alternative, more efficient, approach, though: (i)~first find the $k+1$ roots for $\varphi_{\mathbb Q}(\alpha)$ of \eqref{rooteq}, say $\varphi_1,\ldots,\varphi_{k+1}$, (ii)~then per solution $\varphi_i$ find the three solutions (in $\alpha$) of the equation $\varphi_{\mathbb Q}(\alpha)=\varphi_i$, (iii)~after which one selects, from the resulting $3k+3$ solutions, those with a positive real part (where we know that there are precisely $k+1$ of those). 

\noindent
 A second complication lies in the computation of $\delta_{n,k}^{\mathbb Q}$ (for $n=0,\ldots,k-1$), which requires the numerical evaluation of the higher order derivatives 
 \[\xi_{\mathbb Q,k}^{(n)}(\theta,\omega)= \left.\frac{{\rm d}^n}{{\rm d}\alpha^n}\left(\frac{\alpha -\psi_{\mathbb Q}(\beta)}{\varphi_{\mathbb Q}(\alpha)-\beta}\frac{\beta}{\psi_{\mathbb Q}(\beta)}\right)^k\right|_{\alpha=\theta}.\] 
 We have streamlined the evaluation of these quantities relying on two ideas. 
 (i)~In the first place, it is observed that repeated symbolic derivation of ratios leads to large, unmanageable expressions. To remedy this, first write
 \[\xi_{\mathbb Q,k}(\alpha,\omega)\,\big((\varphi_{\mathbb Q}(\alpha)-\beta))\,\psi_{\mathbb Q}(\beta)\big)^k = \big((\alpha -\psi_{\mathbb Q}(\beta))\,\beta\big)^k,\]
 and then differentiate $n$ times using Leibniz' rule. This yields a recursion for the high-order derivatives of $\xi_{\mathbb Q,k}(\alpha,\omega)$; this we did symbolically using SymPy, the open-source Python library for symbolic computation. (ii)~In the second place, inserting $\alpha=\theta$ into this recursion leads to floating point errors due to division by extremely small numbers. To resolve this, we have inserted $\alpha=\theta\pm \varepsilon$ for a small value of $\varepsilon>0$, and have taken the average of the resulting two numbers. 

\end{itemize}
As mentioned, in the experiments we consider the case of exponentially distributed claims with parameter $\mu$. The parameters chosen are $\lambda= 2$, $r=1.2$, $\sigma^2 = 0.02$, $\mu= 2$, such that $\theta^\star = 0.32875.$

\begin{figure}[htb]
\centering
\includegraphics[scale=0.48, angle=270]{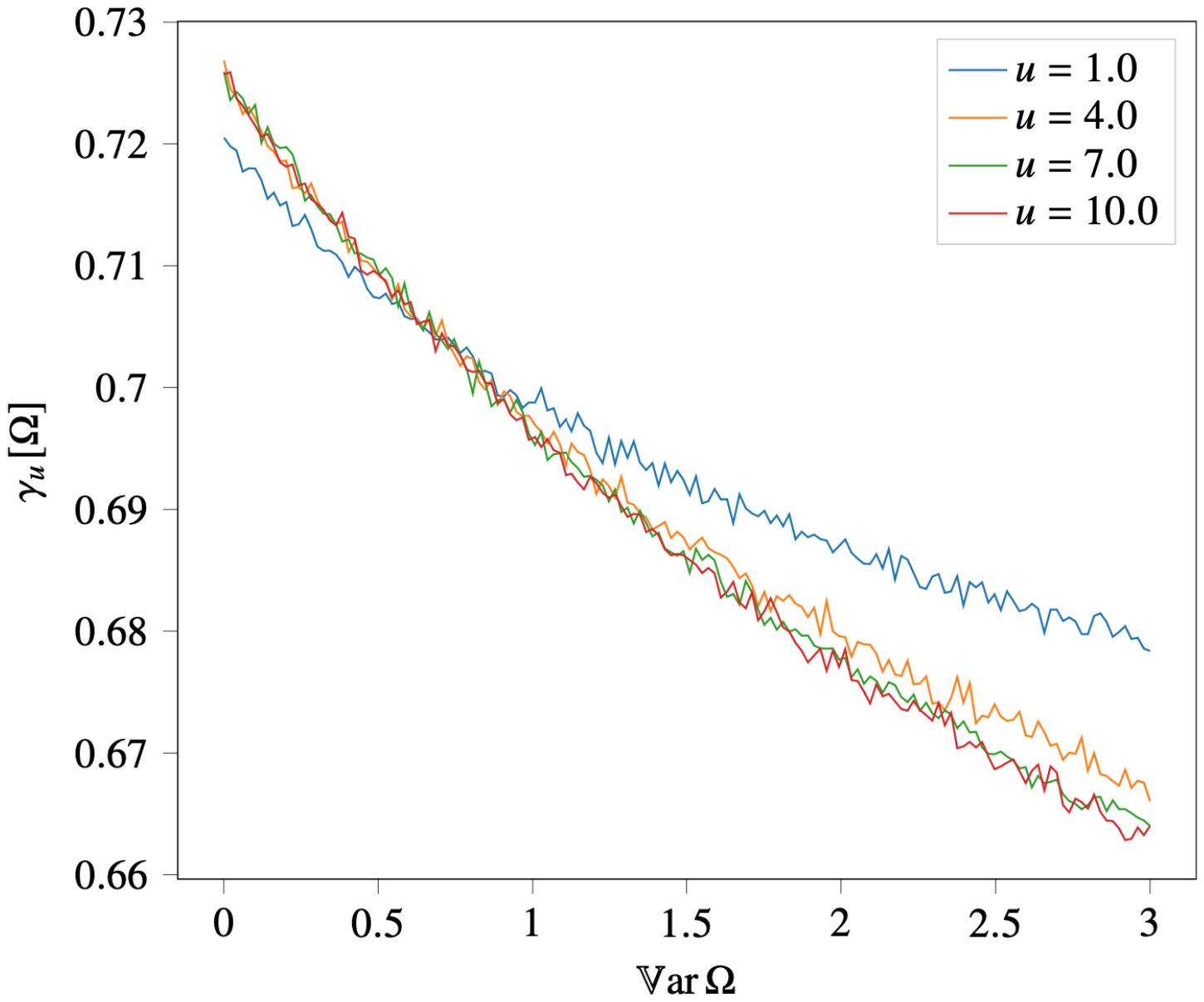}
\caption{Monte Carlo based estimation of $\gamma_u[\Omega]$.}
\label{fig:xxtikz}
\end{figure}

\begin{figure}[htb]
\centering
\includegraphics[scale=0.48, angle=270]{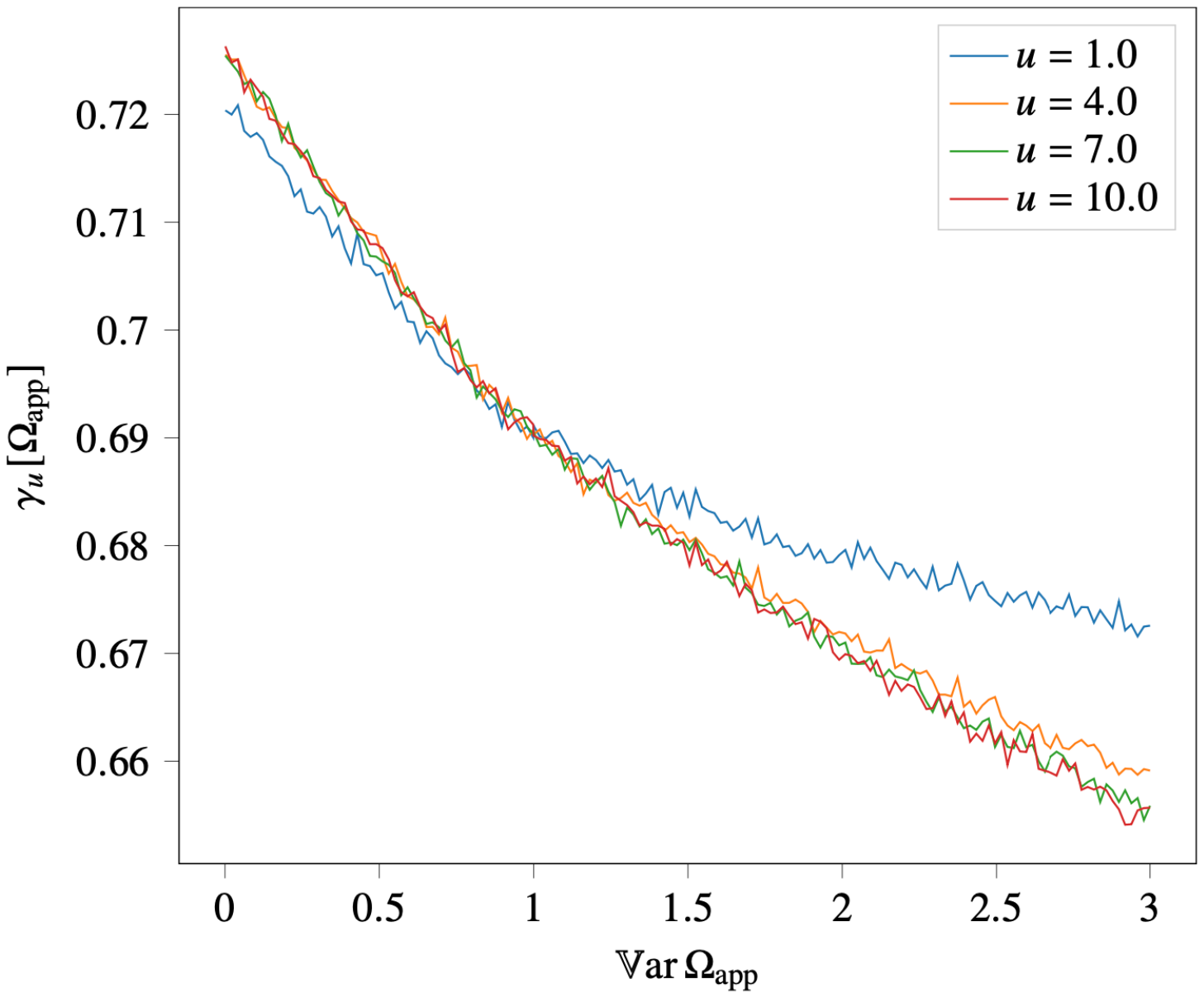}
\caption{Monte Carlo based estimation of $\gamma_u[\Omega_{\rm app}]$.}
\label{fig:yytikz}
\end{figure}

\begin{figure}[htb]
\centering
\includegraphics[scale=0.48, angle=270]{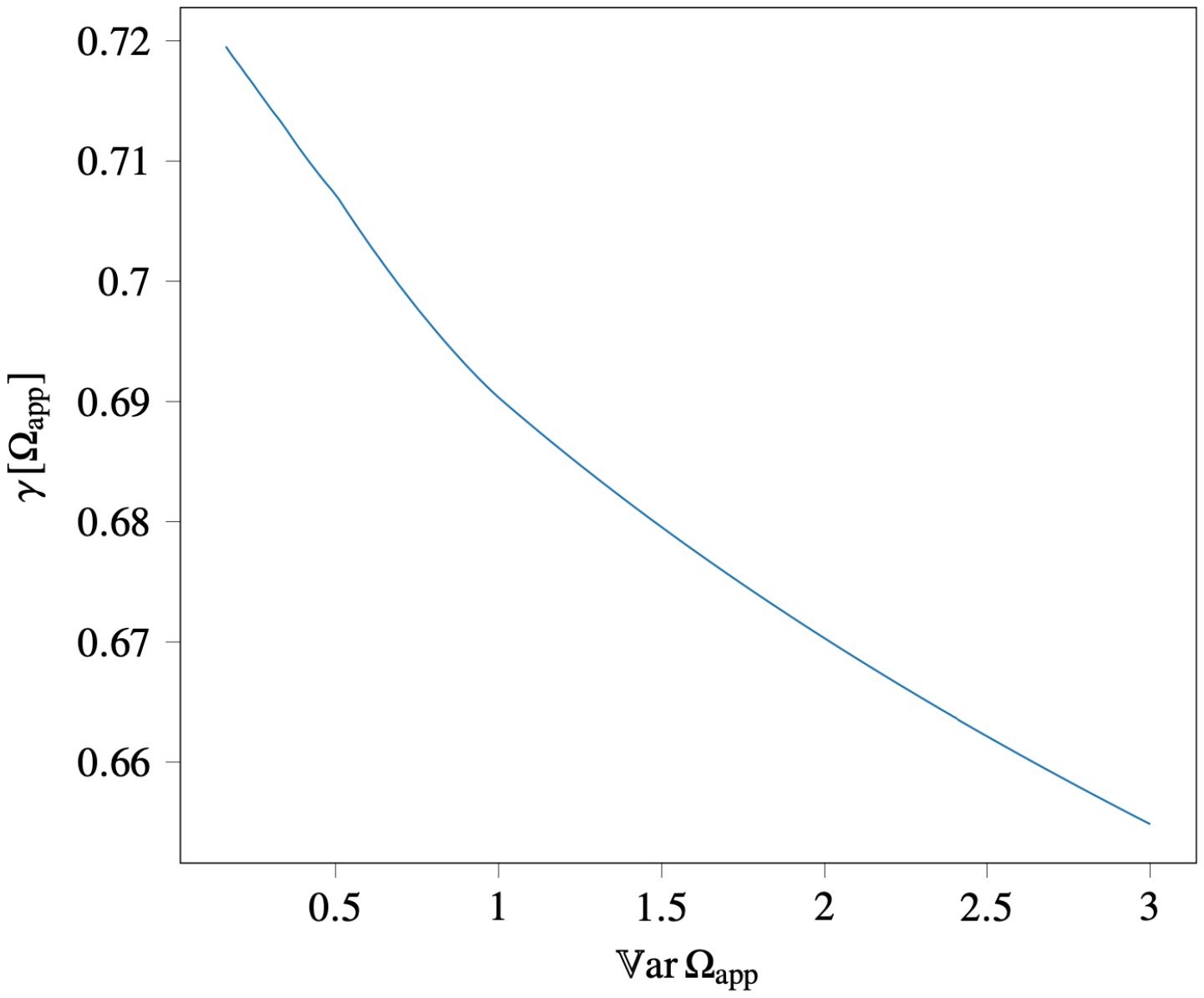}
\caption{Numerical computation for $\gamma[\Omega_{\rm app}]$.}
\label{fig:zztikz}
\end{figure}

In Figure~\ref{fig:xxtikz}  we display, for various values of the initial surplus $u$,  $\gamma_u[\Omega]\coloneqq p(u\,|\,\Omega)\,e^{\theta^\star u}$ as a function of ${\mathbb V}{\rm ar}\,\Omega$, relying on our importance sampling based Monte Carlo procedure. In Figure~\ref{fig:yytikz} we show the corresponding results for its phase-type counterpart $\gamma_u[\Omega_{\rm app}]\coloneqq p(u\,|\,\Omega_{\rm app})\,e^{\theta^\star u}$, again using importance sampling.  Finally, Figure~\ref{fig:zztikz}  displays the curve $\gamma[\Omega_{\rm app}]$, obtained by performing the numerical computations described above. 

All computer code developed for this paper can be accessed in a GitHub repository via the following link: {\small \url{https://github.com/saramorcy/ruin-probabilities}}.
In Figures ~\ref{fig:xxtikz} and~\ref{fig:yytikz} we have displayed, for any value of the variance of the inter-inspection time (${\mathbb V}{\rm ar}\,\Omega$ and ${\mathbb V}{\rm ar}\,\Omega_{\rm app}$, respectively) in the range $\{0.02,0.04,\ldots,2.98,3.00\}$, the quantity $\gamma_u[\Omega]$ resp.\ $\gamma_u[\Omega_{\rm app}]$. Each datapoint is estimated using our importance sampling based Monte Carlo procedure, with the number of runs being equal to $10^5$; we chose this relatively high number of runs to make sure the corresponding standard errors are extremely low, thus making the curves relatively smooth. 

The main conclusions from the figures are the following. 
(i)~The experiments provide empirical evidence that $\gamma_u[\Omega]$ converges to a constant $\gamma[\Omega]$ as $u$ grows large, also for non-phase-type $\Omega$. (ii)~The experiments confirm that $\gamma_u[\Omega_{\rm app}]$ converges to a constant $\gamma[\Omega_{\rm app}]$ as $u$ grows large, in line with our theoretical findings. (iii)~The curves $\gamma[\Omega]$ and $\gamma[\Omega_{\rm app}]$ (that is, as functions of ${\mathbb V}{\rm ar}\,\Omega$ and ${\mathbb V}{\rm ar}\,\Omega_{\rm app}$, respectively)
are hardly distinguishable. This provides backing for the use of the (relatively low dimensional) two-moment approximation of the inter-inspection times. (iv) Already for relatively low values of the initial surplus $u$, the value of $\gamma_u[\Omega_{\rm app}]$ is close to its limiting value $\gamma[\Omega_{\rm app}]$. 

The above observations justify the use of the approximation
\[p(u\,|\,\Omega) \approx \gamma[\Omega_{\rm app}] \,e^{-\theta^\star u}.\]
For the specific parameters chosen in this example, observe that the constant $\gamma[\Omega_{\rm app}]$ only mildly varies in the range ${\mathbb V}{\rm ar}\,\Omega\in[0,3]$, entailing that we could work with  the even simpler approximation
\[p(u\,|\,\Omega) \approx \gamma[T_1] \,e^{-\theta^\star u},\]
where it is recalled that Theorem \ref{thEXP} provides an explicit expression for $\gamma[T_1]$.



\section{Discussion and concluding remarks}\label{DC}
This paper has considered the conventional Cram\'er-Lundberg model, with the distinguishing feature of non-Poissonian inspection. For specific classes of phase-type inter-inspection times the transform of the bankruptcy probability has been identified, where it is noted that by these classes any distribution on the positive half-line can be approximated arbitrarily closely. 
We also establish the tail asymptotics of the bankruptcy probability. By using these as a basis for approximations, which in particular apply in the regime of a high initial surplus, numerical Laplace inversion can be avoided. In addition, a provably optimal  importance sampling algorithm is proposed, applicable to any inter-inspection time distribution. 

Many extensions can be thought of, such as models in which the insurance firm receives interest over its surplus; see e.g.\ \cite[Section VIII.2]{AA} and  \cite{BM} for the case of continuous inspection. Another generalization could concern models with  regime switching; in e.g.\ \cite[Section IX.5]{AA} and \cite{DIKM, DM2} the case of continuous inspection is covered, and \cite{DelMan,KMD} specifically consider the case that the background process is not necessarily irreducible. Finally, one could also consider the model of \cite{CDMR}, in which the environment is periodically resampled. An evident other topic for future research could be the extension to wider classes of distributions for the inter-inspection times, but the discussion in Section \ref{limi} indicates that even for general phase-type distributions there are serious challenges.

\vb

\appendix

\section{Computations for Erlang inter-inspection times}

In order to evaluate the quantity $\pi(\alpha,\beta\,|\,E_{k,\omega})$, we will intensively use the following identities. Some of these are standard, but we include them for completeness. 
\begin{lemma}\label{Erlanglemma} Let $X$ be a nonnegative random variable with LST $\chi(\cdot)$. 

\noindent
(1)~For any $\kappa\in{\mathbb R}$, $w\geqslant 0$ and $k\in{\mathbb N}$,
\[\int_0^w e^{-\kappa v} \frac{\kappa^k v^{k-1}}{(k-1)!}\,{\rm d}v=1 - \sum_{n=0}^{k-1} e^{-\kappa w}\frac{(\kappa w)^n}{n!}.\]

\noindent (2)~For any $\kappa>0$ and $k\in{\mathbb N}_0$,
\[\int_0^\infty e^{-\kappa x} \frac{(\kappa x)^k}{k!}\,{\mathbb P}(X>x)\,{\rm d}x
=\frac{1}{\kappa}\left(1-\sum_{n=0}^k\chi^{(n)}(\kappa)\frac{(-\kappa)^n}{n!}\right),\]
where $\chi^{(n)}(\cdot)$ is the $n$-th derivative of $\chi(\cdot)$.
\end{lemma}
{\it Proof.} Part (1)~follows from the (known) expression for the cumulative distribution function of an Erlang random variable. Part (2) follows by noting that for an Erlang distributed random variable $E_{k,\kappa}$, with shape parameter $k\in\mathbb{N}$ and scale parameter $\kappa$, it holds that
\begin{align*}
{\mathbb P}(X\geqslant E_{k,\kappa}) &= \int_0^\infty {\mathbb P}(E_{k,\kappa}\leqslant x)\,{\mathbb P}(X\in{\rm d}x)\\
&=\int_0^\infty\left(1 - \sum_{n=0}^{k-1} e^{-\kappa x}\frac{(\kappa x)^n}{n!}\right){\mathbb P}(X\in{\rm d}x)=1 -\sum_{n=0}^{k-1}\chi^{(n)}(\kappa)\frac{(-\kappa)^n}{n!},
\end{align*}
where the second equality uses the distribution function of the random variable $E_{k,\kappa}.$ The claim follows by integration by parts.
$\hfill\Box$

\vb

{\it Proof of Lemma \ref{LL2}.}
We start by noticing
\[\int_0^\infty e^{-\alpha u}\,{\mathbb P}(Z_k> u)\,{\rm d}u = \int_0^\infty\int_0^\infty e^{-\alpha u}\frac{\theta^k v^{k-1}}{(k-1)!}e^{-\theta v}{\mathbb P}(Z_k^+>u+v)\,{\rm d}v\,{\rm d}u, \]
where the equality is due to conditioning on the value of $Z_k^-$, which by substituting $w\coloneqq u+v$ and swapping the order of the integrals equals
\begin{equation}\label{f11}\left(\frac{\theta}{\theta-\alpha}\right)^k \int_0^\infty\left( \int_0^w e^{-(\theta-\alpha)v}\frac{(\theta-\alpha)^k v^{k-1}}{(k-1)!}{\rm d}v\right) e^{-\alpha w}\,{\mathbb P}(Z^+_k>w)\,{\rm d}w.\end{equation}
By Lemma \ref{Erlanglemma}.(1), a few elementary steps yield that \eqref{f11} equals
\begin{align*}
    \left(\frac{\theta}{\theta-\alpha}\right)^k\hspace{-0.5mm}\int_0^\infty \hspace{-1mm}e^{-\alpha w}{\mathbb P}(Z_k^+>w)\,{\rm d} w\:-\sum_{n=0}^{k-1}\left(\frac{\theta}{\theta-\alpha}\right)^{k-n}\hspace{-0.5mm}\int_0^\infty \hspace{-1mm}e^{-\theta w}\frac{(\theta w)^n}{n!}\,{\mathbb P}(Z_k^+>w)\,{\rm d}w. 
\end{align*}
Applying  Lemma \ref{Erlanglemma}.(2) to both terms, we find that the integral equals 
$I(\alpha,\beta,\omega)$, as desired. 
$\hfill\Box$

\vb

{\it Proof of Lemma \ref{LL3}.}
We start by splitting the integral under consideration into two terms:  the term ${\pi}_1(\alpha,\beta\,|\,E_{k,\omega})$ corresponding to integration over $\{0\leqslant v\leqslant u\}$ and the term ${\pi}_2(\alpha,\beta\,|\,E_{k,\omega})$ corresponding to  integration over $\{v\leqslant 0\leqslant u\}$; cf.\ the procedure followed in Section \ref{Exp inter-inspection}.

First, we evaluate
\[{\pi}_1(\alpha,\beta\,|\,E_{k,\omega})=\int_0^\infty e^{-\alpha u} \int_{0}^u
{\mathbb P}(Z_k\in {\rm d}v)\,p(u-v,T_\beta\,|\,E_{k,\omega})\,{\rm d}u,\]
by conditioning on the value of $Z_k^-$. We thus obtain
\begin{align*}
&\int_{u=0}^\infty e^{-\alpha u} \int_{v=0}^u \int_{w=0}^\infty
{\mathbb P}(Z_k^+-w\in {\rm d}v)\,e^{-\theta w}\frac{\theta^kw^{k-1}}{(k-1)!}\, p(u-v,T_\beta\,|\,E_{k,\omega})\,{\rm d}w\,{\rm d}u\\
&=\int_{z=0}^\infty\int_{v=0}^\infty e^{-\alpha(v+z)}\int_{w=0}^\infty {\mathbb P}(Z_k^+-w\in {\rm d}v)\,e^{-\theta w}\frac{\theta^kw^{k-1}}{(k-1)!}\,p(z,T_\beta\,|\,E_{k,\omega})\,{\rm d}w\,{\rm d}z
,
\end{align*}
where in the last step the change of variable $z\coloneqq u-v$ has been performed. We see that this expression can be rewritten as
\[\pi(\alpha,\beta\,|\,E_{k,\omega})\int_{v=0}^\infty e^{-\alpha v}\int_{w=0}^\infty {\mathbb P}(Z_k^+-w\in {\rm d}v)\,e^{-\theta w}\frac{\theta^kw^{k-1}}{(k-1)!}\,\,{\rm d}w.\]
Now substituting $v+w\eqqcolon z$ and swapping the order of integration, we obtain
\begin{align*}
\pi&(\alpha,\beta\,|\,E_{k,\omega})\int_0^\infty e^{-\alpha z} \left(\int_0^z e^{-(\theta-\alpha)(z-v)}\frac{\theta^k(z-v)^{k-1}}{(k-1)!}{\rm d}v\right) {\mathbb P}(Z_k^+\in {\rm d}z)\\
&=\pi(\alpha,\beta\,|\,E_{k,\omega})\int_0^\infty e^{-\alpha z} \left(\int_0^z e^{-(\theta-\alpha)v}\frac{\theta^k v^{k-1}}{(k-1)!}{\rm d}v\right) {\mathbb P}(Z_k^+\in {\rm d}z)\\
&=\pi(\alpha,\beta\,|\,E_{k,\omega})\left(\frac{\theta}{\theta-\alpha}\right)^k\int_0^\infty 
e^{-\alpha z} \left(1-\sum_{n=0}^{k-1}e^{-(\theta-\alpha)z}\frac{((\theta-\alpha)z)^n}{n!}\right) {\mathbb P}(Z_k^+\in {\rm d}z)\\
&=\pi(\alpha,\beta\,|\,E_{k,\omega})\left(\frac{\theta}{\theta-\alpha}\right)^k\int_0^\infty 
 \left(e^{-\alpha z}-\sum_{n=0}^{k-1}e^{-\theta z}\frac{((\theta-\alpha)z)^n}{n!}\right) {\mathbb P}(Z_k^+\in {\rm d}z),
\end{align*}
where in the second equality Lemma \ref{Erlanglemma}.(1) has been relied upon. Note that this expression can be interpreted as
\begin{align}\label{ErlangInt2}
\pi(\alpha,\beta\,&|\,E_{k,\omega})\left(\frac{\theta}{\theta-\alpha}\right)^k\left(\xi_k(\alpha,\beta+\omega) - \sum_{n=0}^{k-1} \frac{(\alpha-\theta)^n}{n!} \xi_k^{(n)}(\theta,\beta+\omega)\right)\nonumber\\
\nonumber&=\pi(\alpha,\beta\,|\,E_{k,\omega})\left(\frac{\theta}{\theta-\alpha}\right)^k\left(
\xi_k(\alpha,\beta+\omega) - \sum_{n=0}^{k-1} \left(\frac{\theta-\alpha}{\theta}\right)^n \delta_{n,k}\right)\\
&= \pi(\alpha,\beta\,|\,E_{k,\omega}) \left(\frac{\beta+\omega}{\beta+\omega-\varphi(\alpha)}\right)^k -
\pi(\alpha,\beta\,|\,E_{k,\omega}) \sum_{n=0}^{k-1}\left(\frac{\theta}{\theta-\alpha}\right)^{k-n}\delta_{n,k}.
\end{align}

\begin{figure}[htb]
\centering
\includegraphics[scale=0.58, angle=270]{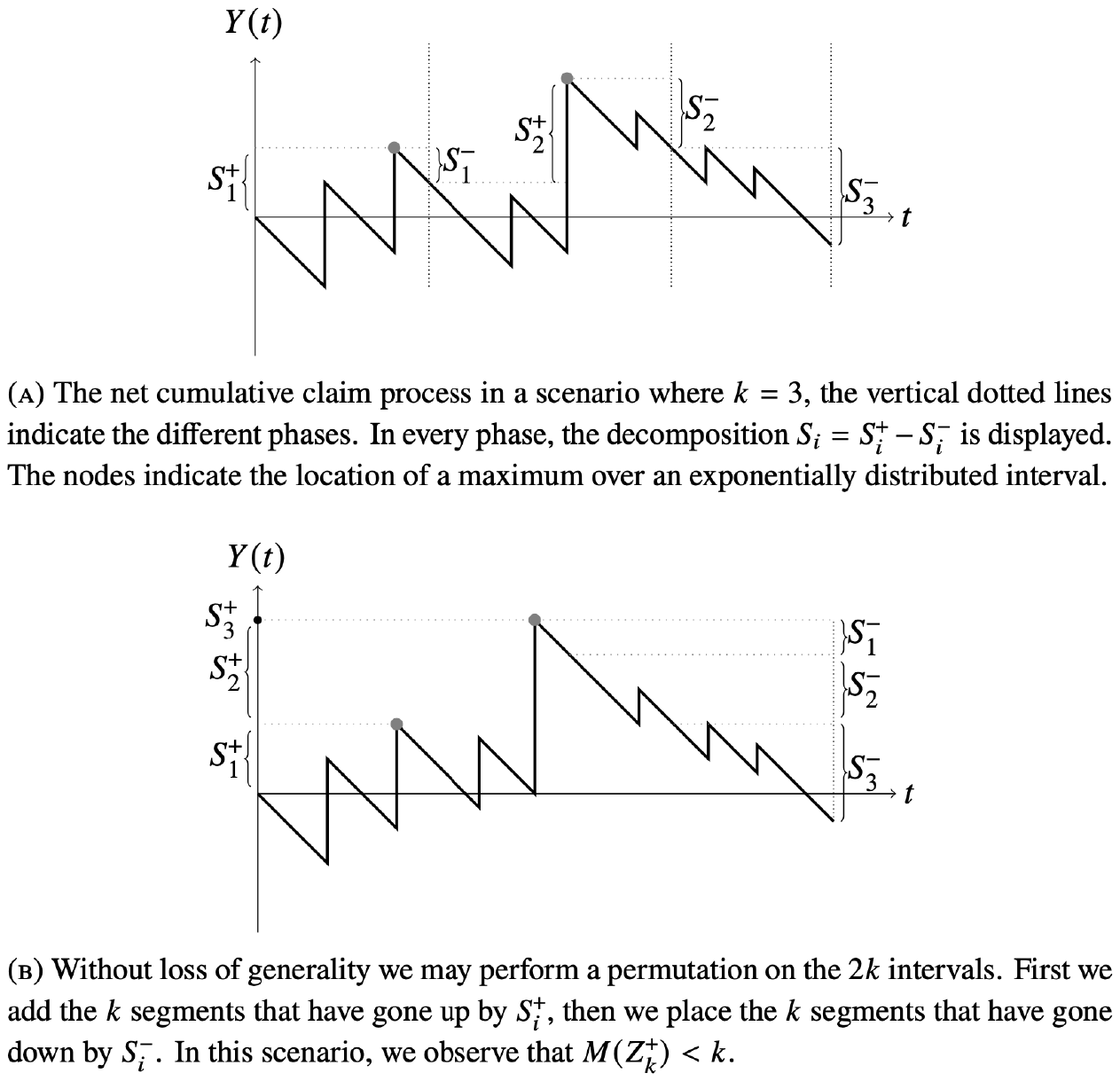}
\caption{The decomposition $Z_k=Z_k^+-Z_k^-$, with $Z_k{=}\sum_{i=1}^k S_i^+ -\sum_{i=1}^k S_i^-$. Depicted is the event $\{Z_k^+<Z_k^-\}$, or equivalently $\{M(Z_k^+)<k\}$.}
\label{->+}
\end{figure}

\begin{figure}[htb]
\centering
\includegraphics[scale=0.58, angle=270]{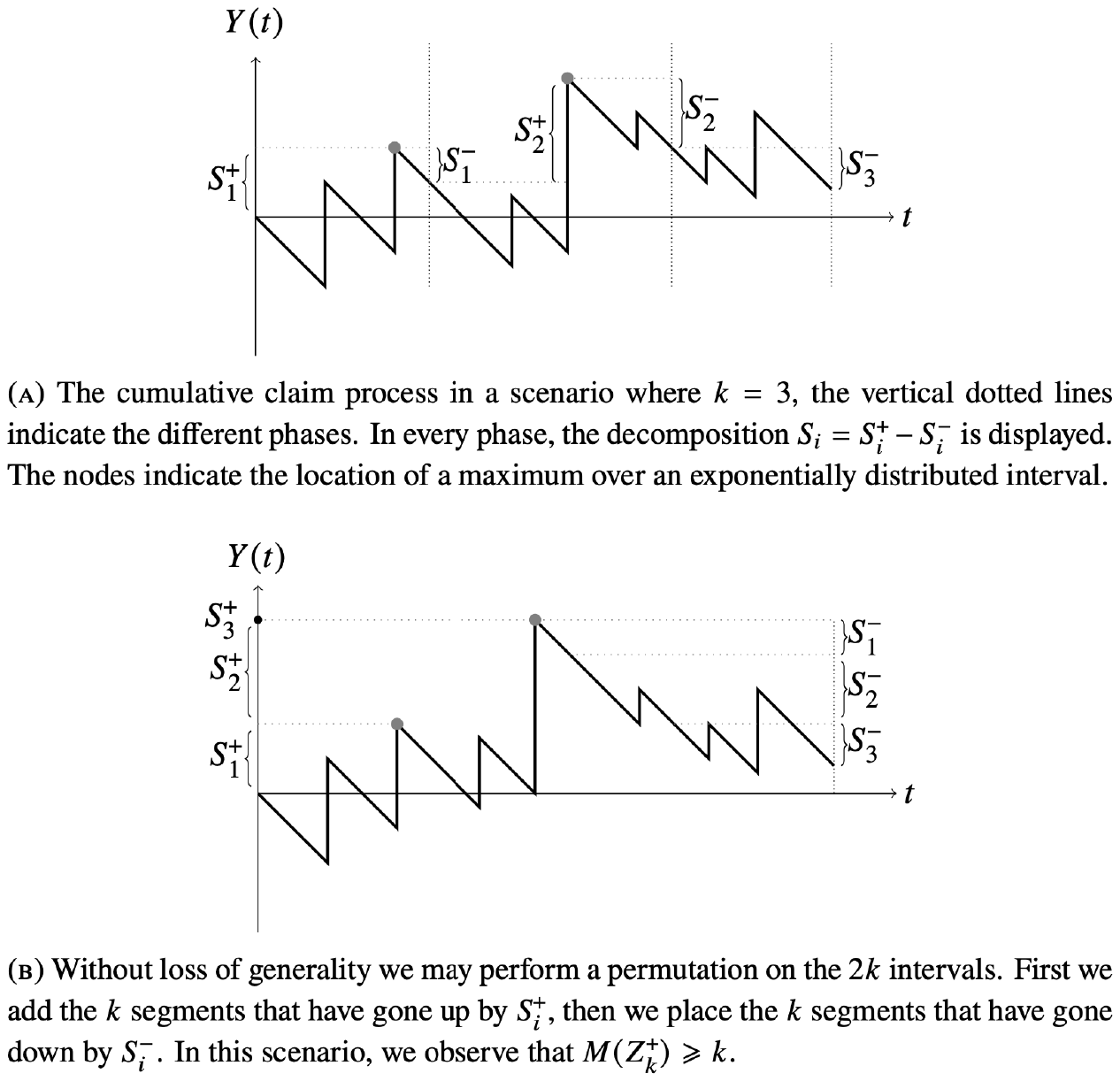}
\caption{The decomposition $Z_k=Z_k^+-Z_k^-$, with $Z_k{=}\sum_{i=1}^k S_i^+ -\sum_{i=1}^k S_i^-$. Depicted is the event $\{Z_k^+>Z_k^-\}$, or equivalently $\{M(Z_k^+)\geqslant k\}$.}
\label{+>-}
\end{figure}

We finally compute
\begin{equation*}\label{I2Erl}{\pi}_2(\alpha,\beta\,|\,E_{k,\omega})=\int_0^\infty e^{-\alpha u} \int_{-\infty}^0
{\mathbb P}(Z_k\in {\rm d}v)\,p(u-v,T_\beta\,|\,E_{k,\omega})\,{\rm d}u.\end{equation*}
Here the first task is to compute ${\mathbb P}(Z_k \in {\rm d} v)$ for $v<0$. To this end, define by $M(t)$ a Poisson process with parameter $\theta$ (independently of $Z_k^+$). Observe that the event $\{Z^-_k>Z^+_k\}$ is equivalent to the event $\{M(Z^+_k)<k\}$, see Figure \ref{->+} and \ref{+>-} (covering the scenarios $M(Z_k^+)<k$ and $M(Z_k^+)\geqslant k$, respectively). Hence, for $v<0$,
\begin{align*}
{\mathbb P}(Z_k \in {\rm d} v)&=
\sum_{\ell=0}^{k-1}{\mathbb P}\big(Z_k \in {\rm d} v, M(Z_k^+)=\ell\big)=\sum_{\ell=0}^{k-1}{\mathbb P}\big(-Z^-_{k-\ell}\in {\rm d} v) \,{\mathbb P}(M(Z_k^+)=\ell)\\
&=\sum_{\ell=0}^{k-1} e^{\theta v}\frac{\theta^{k-\ell}(-v)^{k-\ell-1}}{(k-\ell-1)!}\,{\mathbb P}(M(Z_k^+)=\ell)\,{\rm d}v,
\end{align*}
where the second equality follows by the memoryless property and the last equality follows by realizing that $Z_{k-\ell}^-$ has an Erlang distribution with shape parameter $k-\ell$ and scale parameter $\theta$. Also,
\begin{align*}
{\mathbb P}(M(Z_k^+)=\ell)&=\int_0^\infty {\mathbb P}(M(z)=\ell)\,{\mathbb P}(Z_k^+\in{\rm d}z)=\int_0^\infty e^{-\theta z}\frac{(\theta z)^\ell}{\ell!}\,{\mathbb P}(Z_k^+\in{\rm d}z)\\
&=\frac{(-\theta)^\ell}{\ell!}\xi_k^{(\ell)}(\theta,\beta+\omega)=\delta_{\ell,k}.
\end{align*}

We thus conclude that the distribution of the negative of $Z_k^+-Z_k^-$ can be considered as a mixture of Erlang random variables, entailing that we are to evaluate the integral 
\[\int_0^\infty e^{-\alpha u} \int_{-\infty}^0 \sum_{n=0}^{k-1}e^{\theta v}\frac{\theta^{n+1}(-v)^{n}}{n!}\delta_{k-1-n,k}\,p(u-v,T_\beta\,|\,E_{k,\omega})\,{\rm d}v\,{\rm d}u;\]
note that we swapped the indices by putting $n\coloneqq k-\ell-1$.
Then applying the transformation $w\coloneqq u -v$, and swapping the order of the integrals, after some elementary computations, we arrive at
\begin{align}\label{ErlangInt3}
&\sum_{n=0}^{k-1}\delta_{k-1-n,k}\left(\frac{\theta}{\theta-\alpha}\right)^{n+1} \int_0^\infty\left(\int_0^w e^{-(\theta-\alpha)v}\frac{(\theta-\alpha)^{n+1}v^n}{n!}{\rm d}v\right)\, e^{-\alpha w}\, p(w,T_\beta\,|\,E_{k,\omega})\,{\rm d}w\nonumber \\
&=\sum_{n=0}^{k-1}\delta_{k-1-n,k}\left(\frac{\theta}{\theta-\alpha}\right)^{n+1} \int_0^\infty
 \left(1-\sum_{m=0}^n e^{-(\theta-\alpha)w}\frac{((\theta-\alpha)w)^m}{m!}\right)\,e^{-\alpha w}\, p(w,T_\beta\,|\,E_{k,\omega})\,{\rm d}w\nonumber\\
&=\pi(\alpha,\beta\,|\,E_{k,\omega})\sum_{n=0}^{k-1}\delta_{k-1-n,k}\left(\frac{\theta}{\theta-\alpha}\right)^{n+1} -J(\alpha,\beta,\omega).
\end{align}
Upon adding up the expressions that we have found in \eqref{ErlangInt2} and \eqref{ErlangInt3}, two terms cancel out, and  we obtain the desired expression.
$\hfill\Box$

\end{document}